\documentclass[12pt]{article}
\usepackage{amssymb}
\usepackage{amsmath}

\newtheorem{theo}{Theorem}
\newtheorem{prop}{Proposition}
\newtheorem{lemm}{Lemma}
\newtheorem{coro}{Corollary}
\newtheorem{rema}{Remark}

\newtheorem{conj}{Conjecture}

\newcommand{\cqfd}
{%
\mbox{}%
\nolinebreak%
\hfill%
\rule{2mm}{2mm}%
\medbreak%
\par%
}
\newfont{\gothic}{eufb10}
\date{\empty}
\begin{document}
\title{  Green's  generic syzygy conjecture for curves
of even genus
  lying on a $K3$ surface}
\author{Claire Voisin\\ Institut de math\'ematiques de Jussieu, CNRS,UMR 
7586}

\maketitle
\section{Introduction}
If $C$ is a smooth projective curve of genus $g$ and $K_C$ 
is its canonical bundle,
 the theorem of Noether asserts that the multiplication map
 $$\mu_0: H^0(C,K_C)\otimes H^0(C,K_C)\rightarrow H^0(C,K_C^{\otimes2})$$
 is surjective when $C$ is non hyperelliptic.
 
 The theorem of Petri concerns then the ideal $I$ of $C$ 
 in its canonical embedding, assuming $C$
 is not hyperelliptic. It says that $I$ is generated by its elements of degree
 $2$ if $C$ is neither trigonal nor a plane quintic.
 
 In \cite{Gre}, M. Green introduced and studied the Koszul complexes
 $$\bigwedge^{p+1}H^0(X,L)\otimes H^0(X,L^{q-1})\stackrel{\delta}{\rightarrow}
 \bigwedge^{p}H^0(X,L)\otimes H^0(X,L^{q})\stackrel{\delta}{\rightarrow}
 \bigwedge^{p-1}H^0(X,L)\otimes H^0(X,L^{q+1})$$
 for $X$ a variety and $L$ a line bundle on $X$. Denoting by $K_{p,q}(X,L)$
 the cohomology at the middle of the sequence above, one sees immediately that
 the  surjectivity of the map $\mu_0$ is equivalent to
 $K_{0,2}(C,K_C)=0$, and that if this is the case, the ideal $I$ is
  generated by quadrics
 if and only if $K_{1,2}(C,K_C)=0$. 
 On the other hand, $C$ being non hyperelliptic is equivalent to
 the fact that the {\rm Cliff}ord index ${\rm {\rm Cliff}}\,C$ is strictly positive, where
 $${\rm Cliff} C:=Min\{d-2r,\,\exists L\in Pic\,C, d^0L=d,\,
 h^0(L)=r+1\geq2,\,h^1(L)\geq2\}.$$
 Similarly, $C$ is neither hyperelliptic, nor trigonal nor
  a plane quintic if and only if
 ${\rm Cliff}\,C>1$. 
 
Green's conjecture on syzygies of  canonical curves generalizes then 
the theorems of Noether and Petri as follows
\begin{conj}\label{spe}\cite{Gre} For a smooth projective curve $C$ in 
characteristic $0$,
the condition ${\rm Cliff}\,C>l$ is equivalent to the fact that 
$K_{l',2}(C,K_C)=0,\,\forall l'\leq l$.
\end{conj}
The interest of this formulation of Noether and Petri's theorems is already 
illustrated
in \cite{georoma}, where these theorems are given a modern proof, using
geometric  technics
of computation of 
 syzygies.
 
For our purpose, and as is done in \cite{Gre}, 
it is convenient to use the duality (cf \cite{Gre})
$$K_{p,2}(C,K_C)\cong K_{g-p-2,1}(C,K_C)^*$$
to reformulate the conjecture
as follows
\begin{conj}\label{spedual}\cite{Gre} For a smooth projective curve $C$ 
of genus $g$ in characteristic $0$,
the condition ${\rm Cliff}\,C>l$ is equivalent to the fact that 
$K_{g-l'-2,1}(C,K_C)=0,\,\forall l'\leq l$.
\end{conj}

 If $C$ is now a generic curve, the theorem of Brill-Noether (cf \cite{acgh}, \cite{La}) implies that
 $${\rm Cliff}\,C={\rm gon(C)}-2$$
 where the gonality ${\rm gon(C)}:=Min\,\{d,\,\exists L\in Pic\,C,\,d^0L=d,\,h^0(L)\geq2\}$,
 and that
 $${\rm gon(C)}=\frac{g+3}{2},\,\,if \,\,g\,\,is\,\,odd,$$
 $${\rm gon(C)}=\frac{g+2}{2}
 ,\,\,if \,\,g\,\,is\,\,even.$$
 Hence we arrive at the following conjecture 
 (the generic Green conjecture on syzygies of a canonical curve) :
 \begin{conj}\label{gen} Let $C$ be a generic curve of genus $g$. Then if $g=2k+1$
 or $g=2k$, we have $K_{k,1}(C,K_C)=0$.
 \end{conj}
 \begin{rema} The actual conjecture is $K_{l,1}(C,K_C)=0,\,\forall l\geq k$; 
 but it is easy to prove that 
 $$K_{k,1}(C,K_C)=0\Rightarrow
  K_{l,1}(C,K_C)=0,\,\forall l\geq k.$$
  \end{rema}
 Notice that in the appendix to \cite{Gre}, Green and Lazarsfeld prove the
 conjecture \ref{spe} in the direction $\Leftarrow$ (i. e. they produce 
 non zero syzygies from special linear systems.) Hence the conjecture 
 above cannot be improved,
 namely, under the assumptions above, we have $K_{k-1,1}(C,K_C)\not=0$.
 
 Teixidor \cite{Teix} has recently proposed an approach to the conjecture
 \ref{gen}. Her method uses a degeneration to a tree of elliptic curves and
 the theory of limit linear series of Eisenbud and Harris \cite{EH}, 
 adapted to 
 vector bundles of higher rank.
 It is very likely that her method will lead to a proof of the generic 
 Green conjecture.
 
 We propose here a completely different approach, which does not prove
 Conjecture \ref{gen} in odd genus, but proves Green's Conjecture
 \ref{spedual} for generic curves $C$ of genus $g(C)$ and gonality
 ${\rm gon(C)}$ in the range
 $$\frac{g(C)}{3}+1\leq {\rm gon(C)}\leq\frac{g(C)}{2}+1.$$
 The  inequality on the left says that the gonality has to be not too small,
 but the only  case  which is excluded by the
  inequality on the right is that of a generic curve of odd genus $g$,
 which has gonality $\frac{g+3}{2}$.
  
  Recall from \cite{La} that if $S$ is a $K3$ surface endowed with a ample line bundle
  $L$ such that $L$ generates $Pic\,S$ and $L^2=2g-2$, the smooth members
  $C\in\mid L\mid$ are of genus $g$ and generic in the sense of Brill-Noether, so that 
  in particular they have the same  Clifford
  index as   a generic curve. Hence conjecture \ref{spe} predicts
  that their syzygies vanish as stated in conjecture \ref{gen}.
  This is indeed what we prove here, in the case where the genus is even.
  Note first that the hyperplane restriction theorem \cite{Gre} says that
  \begin{eqnarray}
  \label{0}K_{k,1}(C,K_C)=K_{k,1}(S,L)
  \end{eqnarray}
  whenever $C$ is a hyperplane section of a $K3$ surface $S$ 
  (note that $K_C=L_{\mid C}$
  in this case). Conjecture \ref{gen} for curves of even genus
  is therefore implied by
  \begin{theo}\label{main}
  The pair $(S,L)$ being as above, with $g=2k$, we have
 \begin{eqnarray}
 \label{vanthe}
  K_{k,1}(S,L)=0.
   \end{eqnarray}
  \end{theo}
  
  The body of the paper will be devoted to the
  proof of (\ref{vanthe}).  It turns out that Theorem \ref{main}
  in turn has much stronger
  consequences than the generic syzygy conjecture for curves
  of even genus, and we shall explain this now. In fact we have the following corollary :
  \begin{coro}\label{cor2} For any $\delta\leq \frac{k}{2}$, the generic 
  curve
  of genus 
  $2k-\delta$ which is $k+1-\delta$-gonal satisfies
  $$K_{k,1}(C,K_C)=0$$
  or equivalently by the duality theorem
  $$K_{k-\delta-2,2}(C,K_C)=0.$$
  \end{coro}
  Notice that this result is optimal and exactly predicted by
  Green's conjecture \ref{spe}, since the Clifford index of such
  curve is less than or equal to $k-1-\delta$. An easy computation shows that
  the pairs
  $$g=2k-\delta,\,\rm{gon}=k+1-\delta,\,k\geq0,\,\delta\leq\frac{k}{2}$$
  are exactly the pairs satisfying the inequalities
  $$\frac{g}{3}+1\leq {\rm gon}\leq\frac{g}{2}+1.$$
  Hence Green's Conjecture is proved for generic curves of
  genus and gonality in this range.
  On the other hand, Teixidor \cite{te} has proved the Green conjecture for generic
  curves of fixed gonality in the range
  ${\rm gon}\leq\frac{g}{3}$. Combining this and our corollary
  we see that Green's conjecture is proved for generic curves
  of any fixed gonality, with the exception of generic curves of odd genus
  (which satisfy ${\rm gon}=\frac{g}{2}+1+\frac{1}{2}$).
  
  \vspace{0,5cm}
  
  {\bf  Proof of Corollary \ref{cor2}.} Let
  $(S,L)$ be as in theorem \ref{main}. A generic member
  $X$ of $\mid L\mid$ is $k+1$-gonal. As in section \ref{sec1},
  and following \cite{La},  it follows that there
  is a rank $2$ vector bundle
  $E$ on $S$ with $det\,E=L$, $c_2(E)=k+1$, and $h^0(E)=k+2$.
  The zero set of a generic section of $E$ is a generic member
  of a $g_{k+1}^1$ of a generic curve $X\in \mid L\mid$.
  
  Now let $x_1,\ldots,x_\delta$ be generic points of $S$.
  Beacause $\delta\leq\frac{k}{2}$, the space 
  $$H_{x_\cdot}=H^0(S,E\otimes{\mathcal I}_{x_1}
  \otimes\ldots\otimes{\mathcal I}_{x_\delta})$$
  has rank at least $2$. One checks that 
  for $\alpha,\,\beta$ generic in this space, the curve
  $X$
  defined by the equation
  $$det\,(\alpha\wedge\beta)\in H^0(S,det\,E)=H^0(S,L)$$
  is nodal with nodes exactly as the $x_i$'s.
  On the other hand, the two sections
  $\alpha,\,\beta$ generate a rank $1$ subsheaf
  of the restriction $E_{\mid X}$.
  Let now
  $$n:C\rightarrow X$$
  be the normalization.  The rank $1$ subsheaf
  introduced above induces a line subbundle
  $$D\subset n^*E$$
  with two sections, and it is obvious that the moving part
  of this linear system on $C$ is of degree 
  $k+1-\delta$, since the sections
  $\lambda\alpha+\mu\beta$ of $E$ vanish at the $x_i$'s, so that
  the moving part of their zero sets is of degree
  $k+1-\delta$.
  Hence $C$ is $k+1-\delta$-gonal.
  It remains to show that
  \begin{eqnarray}
  \label{friture}K_{k,1}(C,K_C)=0.
  \end{eqnarray}
  Now we have by the hyperplane restriction theorem and by
  theorem \ref{main} the vanishing
  \begin{eqnarray}
  \label{preder}
 K_{k,1}(X,K_X)=0.
   \end{eqnarray}
   We prove now that this implies (\ref{friture}).
  Notice that there is an identification of 
  $H^0(C,K_C)$ with a subspace of $H^0(X,K_X)$, namely the last space
  is a space of meromorphic forms on $C$ with logarithmic singularities
  over the nodes  satisfying the condition that the sum of the residues
  over each node  vanishes. $H^0(C,K_C)$ is then the subspace of forms which
  are regular.
  
  From this inclusion $j:H^0(C,K_C)\hookrightarrow\ H^0(X,K_X)$, we can deduce a commutative diagram
  of Koszul complexes
  $$
  \begin{matrix}&\bigwedge^{k+1}H^0(C,K_C)&\rightarrow& \bigwedge^kH^0(C,K_C)\otimes H^0(C,K_C)
  &\rightarrow &\bigwedge^{k-1}H^0(C,K_C)\otimes H^0(C,K_C^{\otimes2})\\
  &j\downarrow&&j\downarrow&&j\downarrow\\
  &\bigwedge^{k+1}H^0(X,K_X)&\rightarrow& \bigwedge^kH^0(X,K_X)\otimes H^0(X,K_X)&
  \rightarrow& \bigwedge^{k-1}H^0(X,K_X)\otimes H^0(X,K_X^{\otimes2}).
  \end{matrix}
  $$
  We claim that this induces an inclusion
  $$j:K_{k,1}(C,K_C)\hookrightarrow K_{k,1}(X,K_X).$$
 Indeed,
   consider in general the Koszul
  differential
  $$\delta:\bigwedge^lH^0(Y,{\mathcal L})\rightarrow 
  H^0(Y,{\mathcal L})\otimes\bigwedge^{l-1}H^0(Y,{\mathcal L}).$$
  Then if
  $$\wedge:H^0(Y,{\mathcal L})\otimes\bigwedge^{l-1}H^0(Y,{\mathcal L})
  \rightarrow \bigwedge^lH^0(Y,{\mathcal L})$$
  is the wedge product map,
  one has
  \begin{eqnarray}
  \label{houf}
 \wedge\circ\delta =\pm l Id.
   \end{eqnarray}
  Consider now the inclusion
  $$j:H^0(C,K_C)\otimes\bigwedge^kH^0(C,K_C)
  \rightarrow H^0(X,K_X)\otimes\bigwedge^kH^0(X,K_X).$$
  Let $\alpha\in H^0(C,K_C)\otimes\bigwedge^kH^0(C,K_C)$
  such that $\delta\alpha=0$ and
  $j(\alpha)=\delta\beta$. Then
  (\ref{houf}) gives
  $$j(\alpha)=\delta\beta=\pm\frac{1}{k+1}\delta(\wedge\circ\delta\beta) $$
  $$=\pm\frac{1}{k+1}\delta(\wedge(j(\alpha))).$$
  But $\wedge(j(\alpha))=j(\wedge\alpha)$, so
  that this implies by injectivity
  of $j$ that $\alpha=\pm\frac{1}{k+1}\delta(\wedge\alpha)$. Hence
  $\alpha$ is in fact exact. Hence our claim is proven.
  \cqfd
  In the missed case of a generic curve of odd genus, we have the following
  corollary :
  \begin{coro}\label{cor1} Let $C$ be a generic curve of genus $g=2k-1$; then
  $$K_{k,1}(C,K_C)=0.$$
  \end{coro}
  (Notice that the generic Green conjecture predicts in fact that
  $K_{k-1,1}(C,K_C)=0$.)
  
  \vspace{0,5cm}
  
  {\bf Proof of Corollary \ref{cor1}.} The $K3$ surface
  $S$ being as above, let $X$ be a member of $\mid L\mid$
  with exactly one node as singularity. Let $C$ be the normalization of
  $X$. Then the genus of $C$ is equal to $2k-1$.
  
  We have as before an inclusion
 \begin{eqnarray}
  \label{clouc}
  H^0(C,K_C)\hookrightarrow H^0(X,K_X)
  \end{eqnarray}
  which induces by the same argument as in the proof of corollary
  \ref{cor2} an inclusion
  $$K_{k,1}(C,K_C)\hookrightarrow K_{k,1}(X,K_X).$$
  The hyperplane restriction theorem can be applied to $X\subset S$,
  and together with the vanishing (\ref{vanthe}), it gives
  $$K_{k,1}(X,K_X)=0.$$
  Hence $K_{k,1}(C,K_C)=0$.
  \cqfd
  
 We conclude this introduction
   with a sketch of the main ideas in the proof
   of  theorem \ref{main}. The very starting point 
   is the following observation :
   denote by $S^{[l]}$ the Hilbert scheme parametrizing
   $0$-dimensional length $l$ subschemes of $S$.
   Let $I_l\subset S\times S^{[l]}$ be the incidence subscheme and
  $$\begin{matrix}&I_l&\stackrel{\pi_l}{\rightarrow}&S^{[l]}&\\
  &q\downarrow&&&\\
  &S&&&
  \end{matrix}
  $$
  be the incidence correspondence. Let
  $${\mathcal E}_L:=R^0{\pi_l}_*q^*L$$
  and $L_l:=det\,{\mathcal E}_L$. Then we have
  
  {\bf Fact.} {\it $K_{l-1,1}(S,L)=0$ if and only if
  $$H^0(I_l,\pi_l^*L_l)=\pi_l^*H^0(S^{[l]},L_l).$$}

  Our strategy will be then to construct a subvariety $Z$ of $S^{[k+1]}$, such that
  $$H^0(\tilde Z,\pi_l^*L_l)=\pi_l^*H^0(Z,L_l)$$
 where $\tilde Z:=\pi_l^{-1}(Z)$, and the restriction map
 $$H^0(I_l,\pi_l^*L_l)\rightarrow H^0(\tilde Z,\pi_l^*L_l)$$
 is injective.
 
 As in the papers \cite{La}, \cite{GRH}, the key role 
 in constructing our variety $Z$ and verifying
 the conditions above will be played by the Lazarsfeld-Mukai vector bundle on $S$
 associated with minimal degree base-point free linear systems on smooth members of $\mid L\mid$.

 {\it Terminology.} In this paper, we shall say that a Zariski open subset
 $U\subset X$ is
 large if the complementary closed subset $Z=X-U$ has codimension
 non smaller that $2$ in $X$. In the considered cases, the variety $X$ will be normal, 
 and we will use freely the fact that for a line bundle ${\mathcal L}$
 on $X$
 $$H^0(X,{\mathcal L})\cong H^0(U,{\mathcal L}_{\mid U})$$
 for $U$ a large open subset of $X$.

\section{Strategy of the proof \label{sec1}}
We start with the following observation : Let $X$ be a smooth projective variety. Denote by
$X^{[k]}_{curv}$ the Hilbert scheme parametrizing
curvilinear $0$-dimensional subschemes of $X$ of length $k$.
$X^{[k]}_{curv}$ is smooth, and if $X$ is a curve or a surface, it is a
large open set in the Hilbert scheme $X^{[k]}$ which is smooth.

Let 
$$\begin{matrix}& I_k&\stackrel{\pi_k}{\rightarrow }&X^{[k]}_{curv}\\
&q\downarrow&&\\
&X&&
\end{matrix}
$$
be the incidence correspondence. For a line bundle
$L$ on $X$ denote by ${\mathcal E}_L$ the vector bundle on $X^{[k]}_{curv}$
defined by ${\mathcal E}_L=R^0{\pi_k}_*q^*L$, and let
$$L_k:=det\,{\mathcal E}_L.$$
We have
\begin{lemm} \label{obs} There is a natural isomorphism
$$K_{k,1}(X,L)\cong
H^0(I_{k+1},\pi_{k+1}^*L_{k+1})/\pi_{k+1}^*H^0(X^{[k+1]}_{curv},L_{k+1}).$$ 
In particular, 
$K_{k,1}(X,L)=0$ is equivalent to
$$H^0(I_{k+1},\pi_{k+1}^*L_{k+1})=\pi_{k+1}^*H^0(X^{[k+1]}_{curv},L_{k+1}).$$
\end{lemm}
{\bf Proof.} Recall that $K_{k,1}(X,L)$ is the cohomology at
the middle of the sequence
\begin{eqnarray}
\label{seq}
 \bigwedge^{k+1}H^0(X,L)
\stackrel{\delta}{\rightarrow }H^0(X,L)\otimes \bigwedge^{k}H^0(X,L)
\stackrel{\delta}{\rightarrow }H^0(X,L^{\otimes2})\otimes 
\bigwedge^{k-1}H^0(X,L).
\end{eqnarray}
Now note that there is a natural morphism
\begin{eqnarray}
\label{tau}
\tau: I_{k+1}\rightarrow X\times X^{[k]}_{curv}
\end{eqnarray}
which to $(x,z),\,x\in Supp\, z$ associates $(x,z')$, where
$z'$ is the residual scheme of $x$ in $z$. This morphism is 
well defined because we are working with curvilinear schemes. 

One shows easily that $\tau$ identifies $I_{k+1}$ to a large open subset of
the blow-up of $X\times X^{[k]}_{curv}$ along the
incidence subscheme $I_k$. (Indeed, away from $I_k$ the inverse 
$\tau^{-1}$ of $\tau$ is
given by
$$\tau^{-1}((x,z))=x\cup z,\,x\in X,\,z\in X^{[k]}_{curv}.)$$
Furthermore, if $D\subset I_{k+1}$ is the
 exceptional divisor
one has
\begin{eqnarray}
\label{pulb}
\pi_{k+1}^*L_{k+1}=\tau^*(L\boxtimes L_k)(-D).
\end{eqnarray}
This isomorphism is obtained by studying the natural morphism of vector bundles
over $I_{k+1}$
$$\pi_{k+1}^*{\mathcal E}_{L}\rightarrow \tau^*(pr_1^* E\oplus pr_2^*{\mathcal E}_L),$$
which at a point $\tilde{z}$ such that
$\tau(\tilde{z})=(x,z')$ identifies to the restriction map
$$H^0(L_{\mid\tilde{z}})\rightarrow L_{\mid x}\oplus H^0(L_{\mid z'}).$$
It is immediate to see that the cokernel of this morphism is supported on
the exceptional divisor $D$ of $\tau$ and is of rank $1$ on $D$.

It follows that

\begin{eqnarray}\label{sajan}
H^0(I_{k+1},\pi_{k+1}^*L_{k+1})\nonumber
\\
= {\rm Ker}\,(H^0(X,L)\otimes H^0(X^{[k]}_{curv},L_k)
 \stackrel{rest}{\rightarrow} H^{0}(I_k,{L\boxtimes L_k}_{\mid I_k})).
 \end{eqnarray}

 We now apply the description above to $I_k$ : we note that denoting
 by $p_i,\,i=1,\,2,$ the compositions of the
 projections  with  the
 inclusion $I_k\hookrightarrow X\times X^{[k]}_{curv}$, we have
 $$p_2=\pi_k,\,p_1=pr_1\circ\tau,$$
 where
 $$\tau:I_k\rightarrow X\times X^{[k-1]}_{curv}$$
 is defined as in (\ref{tau}). Hence applying
 formula (\ref{pulb}), we get
 $${L\boxtimes L_k}_{\mid I_k}=\tau^*(L^2\boxtimes L_{k-1})(-D),$$
 where $D$ is now the exceptional divisor of
 the blowing-down morphism
 $\tau:I_k\rightarrow X\times X^{[k-1]}_{curv}$.
So  we conclude  that
 there is a natural inclusion
 $$i:H^0(I_{k},{L\boxtimes L_k}_{\mid I_k})\subset H^0(X,L^{\otimes2})
 \otimes H^0(X^{[k-1]}_{curv},L_{k-1}).$$
 Hence we have constructed from formula (\ref{sajan}) an exact sequence
 $$0\rightarrow H^0(I_{k+1},\pi_{k+1}^*L_{k+1})
 \stackrel{j}{\rightarrow }H^0(X,L)\otimes H^0(X^{[k]}_{curv},L_k)$$
$$ \stackrel{i\,\circ\, rest}{\rightarrow }H^0(X,L^{\otimes2})\otimes 
 H^0(X^{[k-1]}_{curv},L_{k-1}).$$
 Next, it is a standard fact that the natural map
 \begin{eqnarray}
 \label{isowedge}
 \bigwedge^l H^0(X,L)\rightarrow H^0(X^{[l]}_{curv},L_l)
 \end{eqnarray}
 induced by the evaluation map
 $$H^0(X,L)\otimes{\mathcal O}_{X^{[k]}_{curv}}\rightarrow {\mathcal E}_L$$
 are  isomorphisms for any $l$. To check this, one considers the large open
 subset $U^{[l]}$ of $X^{[l]}_{curv}$, made of subschemes
 $z$ which have at most one point of multiplicity $2$ as singularity.
 This set $U^{[l]}$ has the following description : one considers inside
 $X^l$ the large open set $X^l_0$ made of
 $l$-uples $(x_1,\ldots,x_l)$, where at most two of the $x_i$'s coincide. Inside
 $X^l_0$, one blows-up the generalized diagonal
 $\cup_{i\not=j}\{x_i=x_j\}$, which is smooth there,
 and one takes the quotient of the resulting variety $\widetilde{X^l_0}$
 by the action of the symmetric group $S_l$.  By the same argument
 as above one finds that if
 $$r:\widetilde{X^l_0}\rightarrow U^{[l]}$$
 is the quotient map, one has
 $$r^*L_l\cong L\boxtimes\ldots\boxtimes L(-E),$$
 where $E$ is the exceptional divisor of the blowing-up map 
 $\tau':\widetilde{X^l_0}\rightarrow 
 X^l_0$.
 It follows that one has an identification
 $$H^0(X^{[l]}_{curv},L_l)=H^0(U^{[l]},L_l)\cong
 H^0(\tilde{X_0^l},{\tau'}^*(\boxtimes_{i=1}^{i=l}L)(-E))^{S_l}.$$
 The last space is a subspace of
 $$H^0(X^l_0,\boxtimes_{i=1}^{i=l}L)^{S_l}\cong H^0(X^l,\boxtimes_{i=1}^{i=l}L)^{S_l}
 =(\otimes_{i=1}^{i=l}H^0(X,L))^{S_l}.$$
 But looking more closely at the action of
 $S_l$ on the line bundle
 $r^*L_k$, one verifies that the induced action
 of $S_l$ on $\otimes_{i=1}^{i=l}H^0(X,L)$ is the twisted action. 
 Hence the invariant subspace
 $$(\otimes_{i=1}^{i=l}H^0(X,L))^{S_l}$$
 is isomorphic to
 $\bigwedge^lH^0(X,L)$. One verifies then that the injective  map
 $$H^0(X^{[l]}_{curv},L_l)\rightarrow \bigwedge^lH^0(X,L)$$ 
 so obtained is a left inverse for the map (\ref{isowedge}).
 \cqfd
 Using the isomorphisms (\ref{isowedge}), the exact sequence above becomes
 $$0\rightarrow H^0(I_{k+1},\pi_{k+1}^*L_{k+1})
 \stackrel{j}{\rightarrow }H^0(X,L)\otimes \bigwedge^kH^0(X,L)$$
$$ \stackrel{i\,\circ\, rest}{\rightarrow }H^0(X,L^{\otimes2})\otimes 
 \bigwedge^{k-1}H^0(X,L)
 ,$$
 hence it remains it remains only to show that the maps 
 $j\circ\pi_{k+1}^*$ and $i\circ rest$ identify via the isomorphisms
 (\ref{isowedge}) 
  to the differentials $\delta$ of the sequence (\ref{seq}), which is easy. 
 \cqfd
 We consider now a $K3$ surface $S$ endowed with an ample line bundle $L$
 generating $Pic\,S$ and satisfying
 $$L^2=2g-2,\,g=2k,\, k>1.$$
 
 We now explain our strategy to prove the vanishing
 $$K_{k,1}(S,L)=0.$$
  Assume we have
 a subscheme $T\subset S^{[k+1]}$ such that, if $\tilde T$ denotes
 the subvariety $\pi_{k+1}^{-1}(T)$ of $I_{k+1}$, the
  following conditions are satisfied :
 \begin{enumerate}
 \item \label{i}We have an isomorphism
 $$H^0(\tilde T,\pi_{k+1}^*L_{k+1})=\pi_{k+1}^*H^0(T,L_{k+1}).$$
 \item\label{ii} The restriction map
 $$ H^0(I_{k+1},\pi_{k+1}^*L_{k+1})\rightarrow H^0(\tilde T,\pi_{k+1}^*L_{k+1})$$
 is injective.
 \end{enumerate}
 Then we claim that $K_{k,1}(S,L)=0$.

 Indeed we have the trace maps
 $$tr:H^0(I_{k+1},\pi_{k+1}^*L_{k+1})\rightarrow H^0(S^{[k+1]}_{curv},L_{k+1})$$
 $$tr_T:H^0(\tilde T,\pi_{k+1}^*L_{k+1})\rightarrow H^0(T,L_{k+1})$$
 which commute with the restriction maps and which compose to
 $(k+1)\,Id$ with the pull-back maps. If $\sigma\in H^0(I_{k+1},\pi^*L_{k+1})$,
 there exists by  property \ref{i} a $\beta\in H^0(T,L_{k+1})$ such that
 $$\sigma_{\mid\tilde{T}}=\pi_{k+1}^*\beta.$$
 Then 
 $$\beta=\frac{1}{k+1}tr_T(\sigma_{\mid\tilde{T}})=
 (\frac{1}{k+1}tr\,\sigma)_{\mid T}.$$
 Hence
 the section
 $$\sigma'=\sigma-\pi^*(\frac{1}{k+1}Tr\,\sigma)$$
 vanishes on $\tilde T$, hence it is zero by property
 \ref{ii}. So we have  
 $$H^0(I_{k+1},\pi^*L_{k+1})=\pi^*H^0(S^{[k+1]}_{curv},L_{k+1})$$
 and this proves our claim, using lemma \ref{obs}.
 
 We will have to weaken the assumptions above as follows : 
 Suppose we have a normal
  scheme
 $Z$ together with a morphism
 $$j:Z\rightarrow I_{k+1}$$ such that
 $\pi\circ j$ is generically one to one on its image, which is
 not contained in the branch locus of $\pi_{k+1}$. Suppose also that we have
 a normal scheme $Z'$ together with
 a  proper degree $k$ morphism
 $\pi':Z'\rightarrow Z$ and a morphism $j':Z'\rightarrow I_{k+1}$
 satisfying the conditions that
 $$\pi_{k+1}\circ j'=j\circ\pi'$$ and
 the union $j(Z)\cup j'(Z')$ is equal set theoretically to
 $\pi_{k+1}^{-1}(\pi_{k+1}\circ j(Z))$. Finally assume there are subschemes
 $Z'_1\subset Z',\, Z_1\subset Z$ such that
 $$\pi'_{\mid Z'_1}=:\phi:Z'_1\rightarrow Z_1$$
 is a birational isomorphism and
 $j\circ\phi=j'_{\mid Z'_1}$. 
 
 (Hence roughly speaking, and up to birational maps, 
 $\pi_{k+1}^{-1}(\pi_{k+1}\circ j(Z))$ is the scheme obtained by gluing $Z'$ and $Z$ along
 $Z'_1\cong Z_1$.)
 
 Assume now that they satisfy the following set (H) of hypotheses
 \begin{enumerate}
 \item\label{H1} The map
 $${\pi'}^*:H^0(Z,(\pi_{k+1}\circ j)^*L_{k+1})\rightarrow H^0(Z',(\pi_{k+1}\circ j')^*L_{k+1})$$
 is an isomorphism.
 \item\label{H2} The restriction map
 $$H^0(Z,(\pi_{k+1}\circ j)^*L_{k+1})\rightarrow H^0(Z_1,{(\pi_{k+1}\circ j)^*L_{k+1}}_{\mid Z_1})$$
 is injective.
 \item\label{H3} The restriction map
 $$j^*:H^0(I_{k+1},\pi_{k+1}^*L_{k+1})\rightarrow H^0(Z,(\pi_{k+1}\circ j)^*L_{k+1})$$
 is injective.
 \end{enumerate}
 Then we claim that $K_{k,1}(S,L)=0$.
 
 Indeed by Lemma \ref{obs} we have to show that
 $$H^0(I_{k+1},\pi_{k+1}^*L_{k+1})=\pi_{k+1}^*H^0(S^{[k+1]}_{curv},L_{k+1}).$$
 Now if $\sigma\in H^0(I_{k+1},\pi_{k+1}^*L_{k+1})$, by hypothesis
 H\ref{H1}, ${j'}^*\sigma={{\pi'}}^*\alpha$ for some
 $ \alpha\in H^0(Z,(\pi_{k+1}\circ j)^*L_{k+1})$. 
 We show now that $j^*\sigma=\alpha$. Indeed, by
 property H\ref{H2}, it suffices to show that this is true after restriction to 
  $Z_1$, and since
 $\phi:Z'_1\rightarrow Z_1$ is dominating, it suffices to show
 that
 $$\phi^*(\alpha_{\mid Z_1})= {\phi}^*(j^*\sigma_{\mid Z_1}).$$
 But this follows from $j\circ\phi=j'_{\mid Z'_1}$ and from
 ${j'}^*\sigma={\pi'}^*\alpha$, with $\phi=\pi'_{\mid Z'_1}$.
 
 Finally it follows from the equalities  $\alpha=j^*\sigma$
 and ${j'}^*\sigma={\pi'}^*\alpha$ that
 $$\sigma'=\sigma-\pi^*(\frac{1}{k+1}Tr\,\sigma)$$
 vanishes along $j(Z)\cup j'(Z')$. (Indeed, because we know that $j(Z)$
 is not contained in the ramification locus of
 $\pi_{k+1}$ and $\pi_{k+1}\circ j:Z\rightarrow X^{[k+1]}_{curv}$
 is generically one to one on its image, 
 the map
 $\pi_{k+1}$ restricted to $j(Z)\cup j'(Z')$
 is generically a degree $k+1$ unramified map onto its image. On the other hand,
 the two equalities above say that $\sigma_{\mid j(Z)\cup j'(Z')}
 $ is generically a pull-back of a section of $L_{k+1}$ on this image.
 Hence the previously used trace argument applies to an open set of this image.)
 
 Now this implies that $\sigma'$  vanishes by hypothesis
 H\ref{H3}.
 This concludes the proof of our claim.
 \cqfd
 We conclude this section with the description
  of the schemes $Z,\,Z'$ we will be considering.
  
  Recall from \cite{GRH}, \cite{La}, \cite{Mu1}, that there is a unique stable bundle 
  $E$ of rank $2$ on $S$, (the Lazarsfeld-Mukai vector bundle,) which satisfies
  the following properties:
  
  $$det\,E=L,\,
  c_2(E)=k+1,\,
  h^0(E)=k+2.$$
  Such vector bundle is obtained by choosing a line bundle 
  $D$ on a generic member $C$ of $\mid L\mid$, such that
  $h^0(D)=2$ and ${\rm deg}\,D=k+1$. Such a line bundle exists by Brill-Noether
  theory, and it is generated by global sections
  since $C$ does not carry a $g_k^1$ by Lazarsfeld \cite{La}.
  Then we have a vector bundle
  $F$ on $S$ defined by the exact sequence
  \begin{eqnarray}
  \label{exseq}
 0\rightarrow F \rightarrow H^0(D)\otimes{\mathcal O}_S\rightarrow D
 \rightarrow 0
  \end{eqnarray}
  and $E$ is defined as the dual of $F$.
  The stability of $E$ follows from the fact that $Pic\,S={\mathbb Z}L$
  and $H^0(S,E(-L))=0$.
  The uniqueness of such $E$ follows then from the fact that
  $\chi(E,E')=2$ for any other vector bundle $E'$ with the same numerical 
  properties, so that either $Hom\,(E,E')\not=0$ or
  $Hom\,(E',E)\not=0$. But then by stability, $E=E'$.
  
  The property $h^0(S,E)=k+2$ follows from the sequence dual to (\ref{exseq})
  \begin{eqnarray}
  \label{dualexF}
  0\rightarrow H^0(D)^*\otimes{\mathcal O}_S\rightarrow 
  E\rightarrow K_C-D\rightarrow 0,
  \end{eqnarray}
  and from Riemann-Roch which gives
  $h^0(K_C-D)=k$.
  
  Another way to construct the bundle $E$ is via Serre's construction.
   By Riemann-Roch
  the divisors $D$ of degree $k+1$ on smooth members $C$ of $\mid L\mid$
  which satisfy $h^0(C,D)=2$ are
  exactly the subschemes $z$ of degree $k+1$ on $S$ contained in a
  smooth member $C$ of $\mid L\mid$ and satisfying the condition
  that the restriction map
  $$H^0(S,L)\rightarrow H^0(L_{\mid z})$$
  is not surjective.
  Note that since the curves $C$ are general in the sense of Brill-Noether,
  the corank of this map is exactly $1$ and furthermore
  for any $z'\subsetneq z$ the restriction map
  $$H^0(S,L)\rightarrow H^0(L_{\mid z'})$$
  is  surjective.
  Hence, since $K_S$ is trivial,
   to such $z$ corresponds a vector bundle $E$ together with a section
  $\sigma_z$ vanishing on $z$.
  This $E$ is an extension
  \begin{eqnarray}
  \label{***}
  0\rightarrow {\mathcal O}_S\stackrel{\sigma_z}{\rightarrow }E\rightarrow
   \stackrel{\wedge\sigma_z}{\rightarrow } {\mathcal I}_z(L)\rightarrow 0.
   \end{eqnarray}
   Computing the numerical invariants of this bundle $E$, and arguing as before
   by stability, we see that this bundle is isomorphic to the one constructed
   above. Notice that each $g_{k+1}^1$, $D$ on a smooth
   member $C\in \mid L\mid$ provides
   by (\ref{dualexF}) a rank $2$ subspace
   of sections of $E$, and that the zero sets of these sections
   identify to the members of $\mid D\mid$, as subschemes of $S$.
   
   It follows from the exact sequence (\ref{***}) twisted by $E$ that
   $h^0(S,E\otimes{\mathcal I}_z)=1$ for any $z$ as above.
   Hence the morphism
   $${\mathbb P}(H^0(S,E))\rightarrow S^{[k+1]}$$
   which to $\sigma$ associates its zero set,  is in fact an embedding.
   One sees easily that the
   open set ${\mathbb P}(H^0(S,E))_{curv}$ corresponding to
   curvilinear subschemes is large in ${\mathbb P}(H^0(S,E))$.
   
   Let now $W:=\pi_{k+1}^{-1}({\mathbb P}(H^0(S,E))_{curv})\subset I_{k+1}$. $W$ 
   is easily
   shown to be smooth.
   There is a natural morphism
   $$\psi:W\rightarrow S^{[k]}_{curv}$$
   defined as the restriction of
   $pr_2\circ\tau$ to $W$. This $\psi $ can 
   be shown to be generically of degree
   one on its image.
   
   Consider the blow-up $\widetilde{S\times W}$ of $S\times W$ along
   $K:=(Id,\psi)^{-1}(I_k)$. It admits a morphism $\widetilde{(Id,\psi)}$ to the blow-up
   of $S\times S^{[k]}_{curv}$ along $I_k$, and the later 
   contains $I_{k+1}$ as a large open set.
   One verifies that 
   $\widetilde{(Id,\psi)}^{-1}(I_{k+1})$ is a large open set
   of $\widetilde{S\times W}$. This will be our scheme $Z$.
    The morphism $j:Z\rightarrow I_{k+1}$ will be simply
    the restriction to $Z$ of $\widetilde{(Id,\psi)}$.
    
    Again one can show (using now the assumption that $k>1$) that
    the morphism
    $\pi_{k+1}\circ j:Z\rightarrow S^{[k+1]}$ is generically
    of degree one on its image.
    
    Next let $\pi'':\tilde W\rightarrow W$ be the degree $k$ cover 
    obtained by completing
     the Cartesian diagram
     $$
     \begin{matrix}
     &\tilde W&\rightarrow &I_k&\\
     &\pi''\downarrow&&\pi_k\downarrow&\\
    &W&\stackrel{\psi}{\rightarrow } & S^{[k]}_{curv}&
    \end{matrix}.
    $$
    Consider the rational map
    $$j':S\times\tilde W----> I_{k+1}$$
    which to $(s,s_1,w),\,s_1\in Supp\,\psi(w)$ associates
    $(s_1,s\cup \psi(w))$. This morphism
    becomes well defined after blowing-up $K':={(Id,\pi'')}^{-1}(K)$
   and restricting to a large open subset. 
   Our scheme $Z'$ will be this large open set.
     The morphism
     $\pi':Z'\rightarrow Z$ is the restriction to $Z'$
     of the morphism $Bl_{K'}(S\times\tilde W)\rightarrow Bl_K(S\times Z)$ 
     induced by
     $(Id,\pi'')$.
     The morphism $j': Z'\rightarrow I_{k+1}$ is induced by the rational map $j'$ above.
     We have 
     $$\pi_{k+1}\circ j'=\pi_{k+1}\circ j\circ\pi'.$$
     Indeed, both maps send $(s,s_1,w),\,s_1\in Supp\,\psi(w)$ to
     $s\cup\psi(w)$.
     It is obvious that
     $\pi_{k+1}^{-1}(\pi_{k+1}\circ j(Z))$ is equal to $j(Z)\cup j'(Z')$.
     Indeed, the fiber over   
     $s\cup \psi(w)\in\pi_{k+1}\circ j(Z)$ consists in choosing one point in the 
     scheme
     $s\cup \psi(w)$. This point may be $s$, in which case we are in
     $j(Z)$, or has to be a point $s_1$ contained in $Supp\,\psi(w)$ in which case it determines
     a point $(s_1,w)$ of $\tilde W$ over $w$, and we are then in
     $j'(Z')$.
     \begin{rema} The scheme $Z$ is non necessarily smooth, but one can show 
     that
    $K$ is reduced, so that its singular locus is of codimension at least two
    in $S\times W$. The same thing is true for
    $Z'$ and $K'$. If one wants to work with smooth schemes $Z_0$ and $Z'_0$
    (so as to be exactly in the conditions (H) described above), it suffices 
    to restrict
    to the blowing-ups of $S\times W-K_{sing}$ along $K-K_{sing}$
    and $S\times W-K'_{sing}$ along $K'-K'_{sing}$. All what follows
    will be true for these subschemes.
    \end{rema}

     To conclude, it remains now to construct $Z_1$ and $Z'_1$.
     $Z_1$ will be the exceptional divisor of $Z$ (recall that $Z$ is
      a large open set in
     $Bl_K(S\times W)$). Hence $Z_1$ is the inverse image under the
     blow-up map $Z\rightarrow S\times W$ of
     $$K=\{(s,w)\in S\times W,\,s\in Supp\, \psi(w)\}.$$

     We now construct a generic lifting
     of $Z_1$ in $Z'$, the closure of the image of which will be
     $Z'_1$. By definition of
     $Z'$ as a large open set of
     $Bl_{K'}(S\times\tilde W)$, it suffices to construct a lifting
     of  
     $K$ to a component of $K'$ in $S\times\tilde W$. But if
     $(s,w)\in K$, we have $s\in Supp\, \psi(w)$ so that $(s,w)$
      identifies to an element
   $\tilde w$  of $\tilde{W}$. Our lifting sends simply
     $ (s,w)$ to $(s,\tilde w)$.
     
     It remains finally to see that the morphisms $j'$ and $j\circ\pi'$ agree on $Z'_1$.
     Since $I_{k+1}$ is contained in $S\times S^{[k+1]}_{curv}$, it suffices to prove that
     $pr_1\circ j'$ and $pr_1\circ j\circ\pi'$ agree on $Z'_1$ and that
     $pr_2\circ j'$ and $pr_2\circ j\circ\pi'$ agree on $Z'_1$, with
     $pr_2=\pi_{k+1}$ on $I_{k+1}$.
     For the first one this is obvious since
     both maps factor through the contraction
     $Z'_1\rightarrow K'$, and are equal on $K'\subset S\times\tilde W$ 
     to the first projection on $S$, as follows from the definition of the 
     lifting
     $K\rightarrow K'$.
     
     As for the second one, it follows from the fact that, by construction,
     $\pi_{k+1}\circ j'$ and $\pi_{k+1}\circ j\circ \pi'$ agree on $Z'$.
     \cqfd

\section{Proof of the assumptions H\ref{H2} and H\ref{H3}}
We start with the proof of hypothesis H\ref{H2}.
\begin{prop} \label{proviso}Let 
$$Z_1\subset Z\stackrel{\pi_{k+1}\circ j}{\rightarrow }  S^{[k+1]}$$
be as in the previous section. Then the restriction map
$$H^0(Z,(\pi_{k+1}\circ j)^*L_{k+1})\rightarrow
 H^0(Z_1,{(\pi_{k+1}\circ j)^*L_{k+1}}_{\mid Z_1})$$
is injective.
\end{prop}
The proof will be obtained by restricting the construction to a
generic
smooth member $C\in\mid L\mid$. Indeed, recall that
$Z$ is a large open set in the blow-up of $S\times W$ along the
incidence subscheme $K=(id,\psi)^{-1}(I_k)$, where
$$W=\{(x,\sigma)\in S\times {\mathbb P}(H^0(S,E))_{curv},\,\sigma(x)=0\},$$
and $\psi:W\rightarrow S^{[k]}$ sends $(x,\sigma)$ to the residual scheme
of $x$ in $V(\sigma)$. Now since $k\geq1$, the generic element
$z=V(\sigma)$ is supported in a pencil of elements
of $\mid L\mid$, the generic member being smooth. It follows that
a generic element of $S\times W$ is of the form
$(s_1,s_2,z)$, $z=V(\sigma),\,\sigma(s_2)=0$ and there exists a smooth member
$C\in\mid L\mid$ such that $s_1,\,s_2,\,z$ are supported on $C$.
Hence it suffices to prove the analogue of proposition \ref{proviso}
with $Z$ replaced by $Z_C$, the proper transform
of $C\times W_C$ in $Z\subset Bl_K(S\times W)$, where
$$W_C:=\{(c,\sigma)\in C\times {\mathbb P}(H^0(S,E)),\,
\sigma(c)=0,\,V(\sigma)\subset C\},$$
and $Z_1$ is replaced by
$Z_{1,C}:=Z_1\cap Z_C$.

\begin{prop}\label{proh2}
The restriction map

$$H^0(Z_C,{(\pi_{k+1}\circ j)^*L_{k+1}}_{\mid Z_C})
\rightarrow H^0(Z_1,{(\pi_{k+1}\circ j)^*L_{k+1}}_{\mid Z_{1,C}})$$
is injective.
\end{prop}
{\bf Proof.} By the description of the bundle $E$ given in the previous section, we note that
the set 
$$\{\sigma\in{\mathbb P}(H^0(S,E)),
\,\,V(\sigma)\subset C\},$$ identifies by the map $\sigma\mapsto V(\sigma)$
 to the disjoint union 
of the ${\mathbb P}^1\subset C^{(k+1)}$ corresponding to 
$g_{k+1}^1$'s on $C$.  If $D$ is such a $g_{k+1}^1$ on $C$, $D$ gives a morphism of degree
$k+1$
$$\phi_D:C\rightarrow {\mathbb P}^1$$
or a line bundle $L_D$ on $C$ of degree $k+1$ with two sections.
By definition, $W_C$ identifies  (via $\psi$)  to the disjoint union of
copies $C_D$ of $C$ contained in $C^{(k)}$ and indexed
by the $g_{k+1}^1$'s $D$ of $C$, where the map
$$\psi_D:C\cong C_D\rightarrow C^{(k)}$$
is given by 
$$c\mapsto {\rm the\,\,unique\,\, effective\,\,divisor\,\, equivalent
\,\,to}\,\,D-c.$$

Finally $Z_C$ identifies to a disjoint union of surfaces
$Z_{C,D} $ isomorphic to
$C\times C$, since the pull-back $\Delta_D$ to
$C\times C_D$ of  the incidence scheme in $C\times C^{(k)}$ is of 
pure codimension $1$, so that the blow-up of $C\times C_D$
along $\Delta_D$ is isomorphic to $C\times C_D$. Note that under this isomorphism,
the intersection of $Z_1$ with $C\times C_D$ becomes identified to
$\Delta_D$.

Recall now that on the large open set $Z\subset\widetilde{S\times W}$, we have
$$(\pi_{k+1}\circ j)^*L_{k+1}=\tau^*(L\boxtimes\psi^*L_k)(-Z_1).$$
We  have ${L}_{\mid C}=K_C$ and in the sequel
 we will use the notation $H_D$ for the line bundle
${L_k}_{\mid C_D}$. 
(It will be shown that $H_D\equiv kL_D$ but this will not be used now.) Restricting
the equality above to $Z_C$,
we have to show that for each $D$
the restriction map
$$H^0(C\times C, K_C\boxtimes H_D(-\Delta_D))
\rightarrow H^0(\Delta_D,K_C\boxtimes H_D(-\Delta_D)_{\mid\Delta_D})$$
is injective.
In other words we want to show that
\begin{eqnarray}
\label{form28}
H^0(C\times C, K_C\boxtimes H_D(-2\Delta_D))=0.
\end{eqnarray}
Now, since $\Delta_D$ is the restriction
to $C\times C_D$ of the incidence scheme, and since $C_D$ parametrizes
the effective divisors of the form $L_D-x,\,x\in C$, it is clear that
$$\Delta_D=(\phi_D,\phi_D)^{-1}(diag\,({\mathbb P}^1))-diag\,(C).$$
Hence we have
$$\Delta_D\equiv L_D\boxtimes L_D-diag\,(C)$$
in $C\times C$.
It follows that
$$K_C\boxtimes H_D(-2\Delta_D)\equiv
(K_C-2L_D)\boxtimes (H_D-2L_D)+2diag\,C.$$
Now we have the equality
\begin{eqnarray}
\label{forum28}
H^0(C,K_C-2L_D)=0,
\end{eqnarray}
which is proven in \cite{La}, since $C$ is generic in $S$. 
(Indeed for a base point free pencil,
$\mid L_D\mid$, the condition that the $\mu_0$-map
$$H^0(C,L_D)\otimes H^0(C,K_C-L_D)\rightarrow H^0(C,K_C)$$
is injective is equivalent by the base-point free pencil trick to the condition
$$H^0(C,K_C-2L_D)=0.$$
The equality (\ref{form28}) follows now from  (\ref{forum28}) and
 from the fact
that the map $H^0(C,2L_D)\rightarrow H^0(2{L_D}_{\mid 2x})$ is surjective for generic $x$ in $C$.
Hence by Riemann-Roch, $H^0(C,K_C-2L_D)=0$ implies $H^0(C,K_C-2L_D+2x)=0$
for generic $x\in C$. It follows that
$$H^0(C\times C,(K_C-2L_D)\boxtimes (H_D-2L_D)+2diag\,C)=0,$$
which proves the proposition \ref{proh2}, and hence proposition \ref{proviso} is proven.
\cqfd
We turn now to the proof of hypothesis H\ref{H3}.
\begin{prop}\label{proh3} The morphism
$Z\stackrel{j}{\rightarrow} I_{k+1}$ being defined as in the previous
section, the pull-back map
$$j^*:H^0(I_{k+1},\pi_{k+1}^*L_{k+1})\rightarrow
H^0(Z,(\pi_{k+1}\circ j)^*L_{k+1})
$$
is injective.
\end{prop}
The proof  proceeds in several steps, and occupies
the remainder of this section.
Recall that $I_{k+1}$ is a large open set
in the blow-up of
$S\times S^{[k]}$ along the incidence subscheme $I_k$
and that we have the following formula
$$\pi_{k+1}^*L_{k+1}=\tau^*(L\boxtimes L_k)(-D),$$
where $D$  is the exceptional divisor and $\tau$ is the blowing-up map.
Since $Z$ is a large open set in  the proper transform of this blowing-up
under the morphism
$(Id,\psi):S\times W\rightarrow S\times S^{[k]}$,
it suffices to prove
\begin{prop}\label{newform}
The restriction map
$$\psi^*:H^0(S^{[k]},L_k)\rightarrow
H^0(W,\psi^*L_k)$$
is injective.
\end{prop}
In order to prove this proposition, we first show
\begin{lemm} \label{ut}Denoting by
$\pi:W\rightarrow {\mathbb P}(H^0(S,E))$ the restriction of the morphism
$\pi_{k+1}:I_{k+1}\rightarrow S^{[k+1]}$,
we have the formula
$$\psi^*L_k=\pi^*{\mathcal O}_{{\mathbb P}(H^0(S,E))}(k).$$
\end{lemm}
{\bf Proof.} By definition,
$\psi^*L_k=det\,\psi^*{\mathcal E}_{L,k}$, where
the bundle 
${\mathcal E}_{L,k}$ has for fiber $H^0(L_{\mid z})$ at a point
$z\in S^{[k]}$.
Now, if $z\in W$, the scheme $z'=\psi(z)$ has length $k$, hence
the restriction map
$$H^0(S,L)\rightarrow H^0(L_{\mid z'})$$
is surjective. On the other hand
if $z''=\pi(z)$, we have $z'\subset z''$ and the restriction map
$$H^0(S,L)\rightarrow H^0(L_{\mid z''})$$
is not surjective.
Hence we have
$$H^0(S,L\otimes{\mathcal I}_{z'})=H^0(S,L\otimes{\mathcal I}_{z''}),$$
and the fiber of
$\psi^*{\mathcal E}_{L,k}$ at $z$ is canonically isomorphic
to 
$H^0(S,L)/H^0(S,L\otimes{\mathcal I}_{\pi(z)})$.
Hence we have
$$\psi^*L_k=-\pi^*det {\mathcal F},$$
where the bundle ${\mathcal F}$
on ${\mathbb P}(H^0(S,E))$ is the bundle with fiber
$H^0(S,L\otimes{\mathcal I}_{z_\sigma})$ at $\sigma,\,z_\sigma=V(\sigma)$.
Now recall that for each $\sigma$ we have the exact sequence
$$0\rightarrow {\mathcal O}_S\stackrel{\sigma}{\rightarrow}
E\stackrel{\wedge\sigma}{\rightarrow}{\mathcal I}_{z_\sigma}(L)
\rightarrow0.$$
This induces the exact sequence
$$0\rightarrow<\sigma>\rightarrow 
H^0(S,E)
\stackrel{\wedge\sigma}{\rightarrow}H^0(S,{\mathcal I}_{z_\sigma}(L))
\rightarrow0.$$
We conclude immediately from this that 
${\mathcal F}$ fits into the exact sequence
$$0\rightarrow {\mathcal O}_{{\mathbb P}(H^0(S,E))}(-2)
\rightarrow H^0(S,E)\otimes{\mathcal O}_{{\mathbb P}(H^0(S,E))}(-1)
\rightarrow{\mathcal F}\rightarrow0.$$
Since $rank \,\,H^0(S,E)=k+2$, it follows that $det\,{\mathcal F}=
{\mathcal O}_{{\mathbb P}(H^0(S,E))}(-k)$.
\cqfd
It follows from this lemma that we have a natural inclusion
\begin{eqnarray}\label{bas}S^kH^0(S,E)^*\hookrightarrow H^0(W,\psi^*L_k).
\end{eqnarray}
(It will be proven in the next section that this inclusion is in
fact
an isomorphism, but we shall not need this here.)

Our strategy to prove proposition
\ref{newform} will be first to construct an isomorphism
\begin{eqnarray}\label{isoole}
H^0(S^{[k]},L_k)=\wedge^kH^0(S,L)\cong S^kH^0(S,E)^*\end{eqnarray}
and then to show that composed with the inclusion
(\ref{bas}), it is equal, up to a coefficient, to
the pull-back map
$\psi^*$.

\vspace{0,5cm}

{\bf Construction of the isomorphism (\ref{isoole})}. 
We note first that the determinant map
$$det:\bigwedge^2H^0(S,E)\rightarrow H^0(S,det\,E)=H^0(S,L)$$
does not vanish on any element
of rank $2$. Indeed, such an element of rank $2$ is given by a 
subspace $W$ of rank $2$ of $H^0(S,E)$, and if its determinant would vanish
this would imply that $W$ generates a rank $1$ subsheaf of $E$ with at least two sections.
But since $Pic\,S$ is generated by $L$ and
$H^0(S,E(-L))=0$ this is impossible.Hence $det$ provides a
morphism
$$d:G_2\rightarrow {\mathbb P}(H^0(S,L)),$$
where $G_2$ is the Grassmannian of rank two vector subspaces
of $H^0(S,E)$,
or dually
a base-point free linear system
$$K:=H^0(S,L)^*\stackrel{d^*}{\rightarrow}
H^0(G_2,{\mathcal L})=\wedge^2H^0(S,E)^*,$$
where ${\mathcal L}$ is the Pl\"ucker polarization on $G_2$. Notice that
since $rank\,K=2k+1$, and $dim\,G_2=2k$,
$d^*$ has to be injective. Since $K$ is base-point free, we have the exact Koszul complex
on $G_2$
$$0\rightarrow \bigwedge^{2k+1}K\otimes{\mathcal L}^{-(2k+1)}\rightarrow
\ldots\rightarrow K\otimes{\mathcal L}^{-1}\rightarrow{\mathcal O}_{G_2}
\rightarrow 0.$$
We can now tensor this sequence with
$S^k{\mathcal E}$, where the rank $2$ vector bundle ${\mathcal E}$ on
$G_2$ is dual to the tautological rank two subbundle and satisfies
$H^0(G_2,S^k{\mathcal E})=S^kH^0(S,E)^*$.

This provides the exact complex
\begin{eqnarray}\label{Kocon}
0\rightarrow \bigwedge^{2k+1}K\otimes{\mathcal L}^{-(2k+1)}\otimes
S^k{\mathcal E}\rightarrow
\ldots\rightarrow K\otimes{\mathcal L}^{-1}
\otimes S^k{\mathcal E}\rightarrow S^k{\mathcal E}
\rightarrow 0.
\end{eqnarray}
In this complex ${\mathcal K}^.$, the term $S^k{\mathcal E}$ is put in degree
$0$. The hypercohomology ${\mathbb H}^{\,0}(G_2,{\mathcal K}^\cdot)$
vanishes. Now we have a spectral sequence
$$E_1^{p,q}=H^q(G_2,{\mathcal K}^p)\Rightarrow
 {\mathbb H}^{\,p+q}(G_2,{\mathcal
K}^\cdot).$$
It is obvious for degree reasons that all differentials
$d_r$ starting from the term
$E^{0,0}_r$ vanish. 
On the other hand the terms 
$E_1^{p,q}$ with $p+q=-1$ are of the form
$$H^q(G_2,\bigwedge^{q+1}K\otimes{\mathcal L}^{-q-1}\otimes
S^k{\mathcal E}).$$
Using the proposition \ref{appen} proven in the appendix, we see
that these terms are all $0$, except for
$$E_1^{-k-1,k}=H^k(G_2,\bigwedge^{k+1}K\otimes{\mathcal L}^{-k-1}\otimes
S^k{\mathcal E}),$$
which is equal to $\bigwedge^{k+1}K$.
It follows that there is only one non zero differential which 
arrives in some $E_r^{0,0}$, namely
$$d_{k+1}:E_{k+1}^{-k-1,k}\rightarrow E_{k+1}^{0,0}.$$
This implies that
$$E_{k+1}^{0,0}=E_{1}^{0,0}=H^0(G_2,S^k{\mathcal E})=S^kH^0(S,E)^*$$
and that the differential
$d_{k+1}$ above is surjective, since the
 spectral sequence abuts to $0$.
Hence we have build a surjective map
$d_{k+1}$ from a subquotient of
$E_1^{-k-1,k}=\bigwedge^{k+1}K$ to $S^kH^0(S,E)^*$.
Since 
$$dim\,
\bigwedge^{k+1}K=dim\,S^kH^0(S,E)^*=C^{k+1}_{2k+1}$$
 this subquotient must in fact be equal to
$\bigwedge^{k+1}K$ and the map $d_{k+1}$ has to be an isomorphism. 
Finally, since $rank\,K=2k+1$, 
$$\bigwedge^{k+1}K=(\bigwedge^kK)^*=\bigwedge^kH^0(S,L).$$
Hence we have constructed our isomorphism
$$d_{k+1}:\bigwedge^kH^0(S,L)\rightarrow S^kH^0(S,E)^*.$$
\cqfd

To conclude the proof of proposition \ref{newform}, it remains only
to show :
\begin{prop} \label{ilfau}The map $d_{k+1}$ constructed above identifies
up to a coefficient to the map
$$\psi^*:H^0(S^{[k]},L_k)\cong\bigwedge^kH^0(S,L)\rightarrow H^0(W,\psi^*L_k),$$
which takes values in
$ S^kH^0(S,E)^*\subset  H^0(W,\psi^*L_k)$.
\end{prop}
{\bf Proof.} First of all it is clear that
$\psi^*$ takes values in $\pi^*H^0({\mathbb P}(H^0(S,E)),{\mathcal O}(k))
=S^kH^0(S,E)^*$. Indeed, this map is the pull-back map
associated to the morphism
$$W\rightarrow Grass(k+1, H^0(S,L))$$
$$z\mapsto H^0(S,L\otimes{\mathcal I}_{z'}),\,z'=\psi(z).$$
But as mentioned in the proof of lemma \ref{ut}, this morphism factors
through $\pi:W\rightarrow {\mathbb P}(H^0(S,E))$.
More precisely, we  noticed in the proof of lemma \ref{ut}
that the restriction map
$$\psi^*:\wedge^kH^0(S,L)\rightarrow S^kH^0(S,E)^*$$
corresponds to the morphism
$${\mathbb P}(H^0(S,E))\rightarrow Grass(k+1,H^0(S,L))$$
$$\sigma\mapsto det(\sigma\wedge H^0(S,E)).$$
But this morphism is the composition of the morphism
$$\beta:{\mathbb P}(H^0(S,E))\rightarrow Grass(k+1,\bigwedge^2H^0(S,E))$$
$$\sigma\mapsto \sigma\wedge H^0(S,E).$$
and of the rational map 
$$ det:Grass(k+1,\bigwedge^2H^0(S,E))\rightarrow
 Grass(k+1,H^0(S,L))$$
 induced by the determinant map
$det:\bigwedge^2H^0(S,E)\rightarrow H^0(S,L)$.

Next, we note that, with the same spectral sequence argument, and replacing
$K=H^0(S,L)^*\subset \bigwedge^2H^0(S,E)^*$ by the base point 
free linear system
$K'=\bigwedge^2H^0(S,E)^*$ on $G_2$, we could have constructed 
more generally a
surjective map
$$ D_{k+1}:\bigwedge^{k+1}(\bigwedge^2H^0(S,E)^*)
\rightarrow  S^kH^0(S,E)^*,$$
whose restriction to $\bigwedge^{k+1}K$ is equal to
$d_{k+1}$.

Hence  proposition \ref{ilfau} will follow from the following
\begin{lemm} The maps $D_{k+1}$ and
$\beta^*$ from
$\bigwedge^{k+1}(\bigwedge^2H^0(S,E)^*)$ to 
$S^kH^0(S,E)^*$
coincide up to a coefficient.
\end{lemm}
{\bf Proof. } We could argue by $Sl(k+2)$-equivariance. A more direct way
to prove this is to note the following : If $W\subset
\bigwedge^2H^0(S,E)^*$ is a rank $k+1$ vector subspace in general position,
it defines a codimension $k+1$ complete intersection subvariety $G_W$ of $G_2$.
Consider the incidence correspondence
$$\begin{matrix}&P&\stackrel{f}{\rightarrow}& {\mathbb P}(H^0(S,E))\\
&g\downarrow&&\\
&G_2&&
\end{matrix}
$$
Then we have a hypersurface $X_W=f(g^{-1})(G_W)$ of ${\mathbb P}(H^0(S,E))$,
which is easily proven to be of degree $k$.
It is clear that
$$H^0(G_2,S^k{\mathcal E}\otimes{\mathcal I}_{G_W})=
H^0({\mathbb P}(H^0(S,E)),{\mathcal O}_{{\mathbb P}(H^0(S,E))}(k)(-X_W)).$$
On the other hand, from the linear system $W$ we can construct a 
Koszul complex which is a resolution
of ${\mathcal I}_{G_W}$. Hence it is clear that
$$D_{k+1}(\bigwedge^{k+1}W)\subset H^0(G_2,S^k{\mathcal E}\otimes{\mathcal
I}_{G_W}).$$
In other words, if  $\eta$ is a generator of
$\bigwedge^{k+1}W$, $D_{k+1}(\eta)$ is a defining equation of $X_W$ or $0$.
It remains then only to prove that
$\beta^*\eta$ also vanishes on $X_W$. But by definition
$$X_W=\{x\in{\mathbb P}(H^0(S,E)),\,\exists0\not=\gamma\in
{\mathbb P}(H^0(S,E)/<x>),\,x\wedge\gamma\perp W\}.$$
This means that for $x\in X_W$,
the composed map
$$W\hookrightarrow\bigwedge^2H^0(S,E)^*
\rightarrow (x\wedge H^0(S,E))^*$$
is not an isomorphism, hence its determinant vanishes.
But this determinant is equal to
$\beta^*\eta(x)$.
\cqfd

\section{Proof of the assumption H\ref{H1}} 
Recall that we have a Cartesian diagram
$$\begin{matrix}
&Z'&\stackrel{\pi'}{\rightarrow }&Z&\\
&\tau'\downarrow&&\tau\downarrow&\\
&S\times \tilde W&\stackrel{(Id,\pi'')}{\rightarrow }&S\times W&
\end{matrix}
$$
 where the vertical maps $\tau,\,\tau'$ are blow-ups and
 the degree $k$ morphism $\pi''$ fits into the Cartesian diagram
 $$\begin{matrix}
 &\tilde W&\stackrel{}{\rightarrow }&I_k&\\
&\pi''\downarrow&&\pi_k\downarrow&\\
& W&\stackrel{\psi}{\rightarrow }&S^{[k]}_{curv}&
\end{matrix}.
$$
We have the morphisms 
$$j':Z'\rightarrow I_{k+1},\,j:Z\rightarrow I_{k+1}$$ such that
$\pi_{k+1}\circ j'=\pi_{k+1}\circ j\circ \pi'$
and the formula
$$(\pi_{k+1}\circ j)^*L_{k+1}=\tau^*(L\boxtimes \psi^* L_k)(-D)$$
where $D$ is the exceptional divisor of $\tau$. Pulling-back this equality
to $Z'$, we have
$$(\pi_{k+1}\circ j')^*L_{k+1}={\tau'}^*(L\boxtimes (\psi\circ \pi'')^* L_k)(-D'),
$$
where $D'$ is the exceptional divisor of
 $\tau'$.
Since $D'={\pi'}^{-1}(D)$ and $\pi'$ is surjective, we conclude that
in order to prove H\ref{H1}, that is the fact that the
pull-back map 
$${\pi'}^{*}:H^0(Z,(\pi_{k+1}\circ j)^*L_{k+1})\rightarrow H^0(Z',(\pi_{k+1}\circ j')^*L_{k+1})$$
is surjective, it suffices to show that 
the pull-back map
$${\pi''}^*:H^0(W,\psi^*L_k)\rightarrow H^0(\tilde W,(\psi\circ \pi'')^*L_k)$$
is surjective.

Now recall that we have a morphism
$$\pi:W\rightarrow {\mathbb P}(H^0(S,E))$$
such that (cf lemma \ref{ut})
$$\psi^*L_k=\pi^*{\mathcal O}_{{\mathbb P}(H^0(S,E))}(k).$$
Denoting by $r:=\pi\circ \pi'':\tilde W\rightarrow {\mathbb P}((H^0(S,E))$,
 we shall prove the following stronger statement
\begin{theo} \label{thmain}The pull-back map
\begin{eqnarray}
\label{fleche}
r^*:H^0({\mathbb P}(H^0(E)),{\mathcal O}_{{\mathbb P}(H^0(E))}(k))
\rightarrow H^0(\tilde W,(\psi\circ \pi'')^*L_k)
\end{eqnarray}
is surjective.
\end{theo}
The end of this section will be devoted to the proof of this theorem, which
proceeds in several steps. In what follows, we shall use the notation
$H^0(E)$ for $H^0(S,E)$.

Notice to begin with that
$\tilde{W}$ is a large open set in the subscheme
$$W'\subset \widetilde{S\times S}\times {\mathbb P}(H^0(E)),$$
where $\widetilde{S\times S}$ is the blow-up of $S\times S$ along the diagonal,
defined as
$$W':=\{(x,y,\eta,\sigma),\,\sigma_{\mid\eta=0},\,\{x\},\,\{y\}\subset\eta\}.$$
(Here $\eta$ is a subscheme of length $2$ of $S$, and we see elements
of $\widetilde{S\times S}$ as elements $(x,y)$ of $S\times S$ together with a schematic structure
$\eta$
of length $2$
on $\{x\}\cup\{y\}$.)

The map $r$ is just the restriction to $W'$ of the second projection. Hence we have
$$H^0(\tilde W,(\psi\circ \pi'')^*L_k)=H^0(W',pr_2^*{\mathcal O}_{{\mathbb P}(H^0(E))}(k))$$
and the surjectivity of (\ref{fleche}) is equivalent to the condition
\begin{eqnarray}
\label{vanh1}H^1(\widetilde{S\times S}\times {\mathbb P}(H^0(E)),pr_2^*{\mathcal O}(k)\otimes
{\mathcal I}_{W'})=0.
\end{eqnarray}
Now notice that there is a vector bundle
$\tilde E_2$ on $\widetilde{S\times S}$ such that $W'$ is the zero set
of a section
$\sigma$ of $\tilde E_2\boxtimes{\mathcal O}_{{\mathbb P}(H^0(E))}(1)$. Indeed
it suffices to take for $\tilde{E}_2$ the vector bundle with fiber
$H^0(E_{\mid \eta})$ at the point $(x,y,\eta)$ of $\widetilde{S\times S}$.
Then the section $\sigma$ takes the value
$\tau_{\mid\eta}$ at the point $(x,y,\eta,\tau)$ of $\widetilde{S\times S}\times {\mathbb P}(H^0(E))$.
One checks easily that $W'$ is reduced of codimension $4$. Hence we have
a Koszul resolution of ${\mathcal I}_{W'}$
\begin{eqnarray}
\label{kore}
0\rightarrow \bigwedge^4{\tilde E_2}^*\boxtimes{\mathcal O}(-4)\rightarrow \ldots\rightarrow 
{\tilde E_2}^*\boxtimes{\mathcal O}(-1)\rightarrow{\mathcal I}_{W'}
\rightarrow 0.
\end{eqnarray}
Our first goal will be to compute the cohomology groups
of $\widetilde{S\times S}\times {\mathbb P}(H^0(E))$
with value in $\bigwedge^i{\tilde E_2}^*\boxtimes{\mathcal O}(k-i)$. Since $k\geq2,\,i\leq4$,
${\mathcal O}(k-i)$ has no higher cohomology on
${\mathbb P}(H^0(E))={\mathbb P}^{k+1}$. Hence we have
$$H^l(\widetilde{S\times S}\times {\mathbb P}(H^0(E)),
\bigwedge^i{\tilde E_2}^*\boxtimes{\mathcal O}(k-i))=
H^l(\widetilde{S\times S},\bigwedge^i{\tilde E_2}^*)\otimes S^{k-i}H^0(S,E)^*.$$
We have now the following proposition
\begin{prop}\label{calco}\begin{enumerate}
\item\label{iti}
$H^2(\widetilde{S\times S},{\tilde E_2}^*)=pr_1^*H^2(S,E^*)\oplus pr_2^* H^2(S,E^*)$ and
$$ H^1(\widetilde{S\times S},{\tilde E_2}^*)=0.$$
\item\label{itii}
$H^2(\widetilde{S\times S},\bigwedge^2{\tilde E_2}^*)=pr_1^*H^2(S,-L)\oplus pr_2^* H^2(S,-L).$
\item\label{itiii} $H^4(\widetilde{S\times S},\bigwedge^4{\tilde E_2}^*)$
is dual to ${\rm Ker}\,(H^0(S,L)\otimes H^0(S,L)\rightarrow H^0(S,2L))$.
\item\label{itiv} $H^3(\widetilde{S\times S},\bigwedge^3{\tilde E_2}^*)=0$ and
$H^4(\widetilde{S\times S},\bigwedge^3{\tilde E_2}^*)$ admits as a quotient
$$H^4(\widetilde{S\times S},\tau^*(E^*\boxtimes (-L))(2\Delta))
\oplus H^4(\widetilde{S\times S},\tau^*((-L)\boxtimes E^*)(2\Delta)),$$ which is dual to 
the direct sum of two copies of
$$ {\rm Ker}\,(H^0(S,E)\otimes H^0(S,L)\rightarrow H^0(S,E\otimes L)).$$
(Here $\Delta\subset \widetilde{S\times S}$ is the exceptional divisor.)
\end{enumerate}
\end{prop}
{\bf Proof.}

\ref{iti}. The bundle $\tilde E_2$ fits into the exact sequence
\begin{eqnarray}
\label{exE2}
0\rightarrow \tilde E_2\rightarrow \tau^*(pr_1^*E\oplus pr_2^*E)\rightarrow {\tau'}^*E\rightarrow 0
\end{eqnarray}
where $\tau:\widetilde{S\times S}\rightarrow S\times S$ is the blowing-down map, and
where $\tau':\Delta\rightarrow Diag\,S$ is its restriction to the exceptional divisor.

Dualizing, we get the following exact sequence 
\begin{eqnarray}
\label{exE2*}0\rightarrow \tau^*(pr_1^*E^*\oplus pr_2^*E^*)
\rightarrow {\tilde E_2}^*\rightarrow {\tau'}^*E^*\otimes
{\mathcal O}_\Delta(\Delta)\rightarrow 0.
\end{eqnarray}
Now $R^0\tau'_*{\mathcal O}_\Delta(\Delta)=R^1\tau'_*{\mathcal O}_\Delta(\Delta)=0$
hence the sheaf on the right has no cohomology. It follows that
$$H^i(\widetilde{S\times S},{\tilde E_2}^*)=H^i(\widetilde{S\times S},\tau^*(pr_1^*E^*\oplus pr_2^*E))$$
$$=H^i(S\times S,pr_1^*E^*\oplus pr_2^*E).$$
Since $E^*$ has no odd dimensional cohomology, nor ${\mathcal O}_S$, it follows from
K\"unneth formula that the same is true for $pr_1^*E^*\oplus pr_2^*E$ on $S\times S$.
Finally we have
$$H^2(S\times S,pr_1^*E^*)=H^2(S,E^*)$$
since $H^0(S,E^*)=0$. This proves \ref{iti}.

\ref{itii}. From (\ref{exE2*}) we deduce that
$\bigwedge^2\tilde E_2^*$ has a filtration whose successive terms are
$$\bigwedge^2\tau^*(pr_1^*E^*\oplus pr_2^* E^*),\,
\tau^*(pr_1^*E^*\oplus pr_2^*E^*)\otimes{\tau'}^*E^*\otimes{\mathcal O}_\Delta(\Delta),\,
\bigwedge^2{\tau'}^*E^*\otimes{\mathcal O}_\Delta(2\Delta).$$
The sheaf $(pr_1^*E^*\oplus pr_2^*E^*)\otimes{\tau'}^*E^*\otimes{\mathcal O}_\Delta(\Delta)$
has no cohomology, since ${\mathcal O}_\Delta(\Delta)$ has no cohomology along the fibers of
$\tau'$. Hence we have an exact sequence
$$H^1(\Delta,\bigwedge^2{\tau'}^*E^*\otimes{\mathcal O}_\Delta(2\Delta))\rightarrow 
H^2(\widetilde{S\times S},\tau^*\bigwedge^2(pr_1^*E^*\oplus pr_2^* E^*))\rightarrow 
H^2(\widetilde{S\times S},\bigwedge^2\tilde E_2^*)$$
$$
\rightarrow  H^2(\Delta,\bigwedge^2{\tau'}^*E^*\otimes{\mathcal O}_\Delta(2\Delta))
\ldots
 $$
 But since 
 $$R^1\tau'_*(2\Delta_{\mid \Delta})={\mathcal O}_S,\,R^0\tau'_*(2\Delta_{\mid \Delta})=0,$$
 the term on the left is equal to $H^0(S,\bigwedge^2E^*)=0$ and the term on the right is equal to
 $H^1(S,\bigwedge^2E^*)=0$.
 Hence we have
 $$H^2(\widetilde{S\times S},\bigwedge^2\tilde E_2^*)=
 H^2(\widetilde{S\times S},\tau^*\bigwedge^2(pr_1^*E^*\oplus pr_2^* E^*))$$
 $$=
  H^2({S\times S},\bigwedge^2(pr_1^*E^*\oplus pr_2^* E^*))$$
 Finally 
 $$\bigwedge^2(pr_1^*E^*\oplus pr_2^* E^*)=pr_1^*\bigwedge^2E^*\oplus E^*\boxtimes E^*
 \oplus pr_2^*\bigwedge^2E^*.$$
 The central term has no cohomology in degree $2$ by K\"unneth formula, because
 $H^1(S,E^*)=H^0(S,E^*)=0$, and we have
 $$H^2({S\times S},pr_1^*\bigwedge^2E^*)=H^2(S,\bigwedge^2E^*)=H^2(S,-L),$$
 where the first equality follows from K\"unneth formula and
 $H^0(S,\bigwedge^2E^*)=0$.
 This proves \ref{itii}.
 
 \ref{itiii}.
  We have $det\,\tilde E_2^*=\tau^*((-L)\boxtimes(-L))(2\Delta)$ by the exact sequence
 (\ref{exE2}). Hence
 \begin{eqnarray}
 \label{eqsa}
 \bigwedge^3\tilde E_2^*=\tilde E_2\otimes det\,\tilde E_2^*=
 \tilde E_2\otimes\tau^*((-L)\boxtimes(-L))(2\Delta).\end{eqnarray}
 The exact sequence
 (\ref{exE2}) gives now the long exact sequence
 $$ H^2(\Delta,{\tau'}^*(E(-2L))(2\Delta_{\mid\Delta}))
 \rightarrow  H^3(\widetilde{S\times S},\bigwedge^3\tilde E_2^*)$$
 $$
 \rightarrow
  H^3(\widetilde{S\times S},\tau^*((pr_1^*E\oplus pr_2^*E)\otimes((-L)\boxtimes(-L)))(2\Delta)).$$
 Since $R^0\tau'_*{\mathcal O}_\Delta(2\Delta)=0,\,R^1\tau'_*{\mathcal O}_\Delta(2\Delta)=
 {\mathcal O}_S$, the left hand side is equal to
 $H^1(S,E(-2L))$, which is easily seen to be $0$. 
 
 Next we have $K_{\widetilde{S\times S}}={\mathcal O}_{\widetilde{S\times S}}(\Delta)$, hence
 $$H^3(\widetilde{S\times S},\tau^*(pr_1^*E\oplus pr_2^*E)\otimes((-L)\boxtimes(-L)))(2\Delta))$$
 is dual to 
 \begin{eqnarray}
 \label{group}
  H^1(\widetilde{S\times S},\tau^*((pr_1^*E^*\oplus pr_2^*E^*)\otimes(L\boxtimes L))(-\Delta)).
 \end{eqnarray}
 Writing the exact sequence
 $$0\rightarrow \tau^*((pr_1^*E^*\oplus pr_2^*E^*)\otimes(L\boxtimes L))(-\Delta)
 \rightarrow \tau^*((pr_1^*E^*\oplus pr_2^*E^*)\otimes(L\boxtimes L))
 $$
 $$\rightarrow \tau^*((pr_1^*E^*\oplus pr_2^*E^*)\otimes(L\boxtimes L))_{\mid \Delta}
 \rightarrow 0
 $$
 where
 $$\tau^*((pr_1^*E^*\oplus pr_2^*E^*)\otimes(L\boxtimes L))_{\mid \Delta} \cong
 {\tau'}^*((E^*\oplus E^*)(2L)),$$
 and using the isomorphism
 $$E^*\otimes L\cong E,$$
  we see that  the vanishing of the cohomology group
 (\ref{group}) follows from the fact that
  the multiplication
 map
 $$H^0(S,E)\otimes H^0(S,L)\rightarrow H^0(S,E\otimes L)$$
 is surjective, which is easily checked,  and from the vanishing
 $$H^1({S\times S},(pr_1^*E^*\oplus pr_2^*E^*)\otimes L\boxtimes L)=0.$$

 Finally the equality (\ref{eqsa}) and the exact sequence
 (\ref{exE2}) also show that
 $H^4(\widetilde{S\times S},\bigwedge^3\tilde E_2^*)$ admits 
$$ H^4(\widetilde{S\times S},\tau^*((pr_1^*E\oplus pr_2^*E)\otimes((-L)\boxtimes(-L)))(2\Delta))$$
 as a quotient. By Serre's duality this space is dual to
 \begin{eqnarray}
 \label{type}
 H^0(\widetilde{S\times S},\tau^*((pr_1^*E^*\oplus pr_2^*E^*)\otimes(L\boxtimes L))(-\Delta)).
 \end{eqnarray}
 But this is equal to 
 $$H^0({S\times S},(pr_1^*E^*\oplus pr_2^*E^*)\otimes(L \boxtimes L)\otimes{\mathcal I}_{Diag}).$$
 We use then  the fact that
 $$pr_1^*E^*\otimes(L\boxtimes L)=E\boxtimes L$$
 to conclude that 
 (\ref{type}) is equal to the sum of two copies of
 $$  {\rm Ker}\,(H^0(S,E)\otimes H^0(S,L)\rightarrow H^0(S,E\otimes L)).$$
\ref{itiv}.  We already noticed that 
$$\bigwedge^4\tilde E_2^*=det\,\tilde E_2^*=\tau^*((-L)\boxtimes(-L))(2\Delta).$$
It follows then from Serre's duality and 
$K_{\widetilde{S\times S}}={\mathcal O}_{\widetilde{S\times S}}(\Delta)$
that
$H^4(\widetilde{S\times S},\bigwedge^4\tilde E_2^*)$ is dual to
$$H^0(\widetilde{S\times S},\tau^*(L\boxtimes L)(-\Delta))= {\rm Ker}
\,(H^0(S,L)\otimes H^0(S,L)\rightarrow H^0(S,2L))).$$
Hence \ref{itiv} is proven.

\cqfd
Coming back to the Koszul resolution of ${\mathcal I}_{W'}\otimes pr_2^*{\mathcal O}(k)$
induced by 
(\ref{kore}), we see that in order to prove the vanishing
(\ref{vanh1}), it suffices to show :

{\it a) $H^1(\widetilde{S\times S}\times{\mathbb P}(H^0(E)),{\tilde E_2}^*\boxtimes{\mathcal O}(k-1))=0$.

b) The interior product with $\sigma$
$$int(\sigma):H^2(\widetilde{S\times S}\times{\mathbb P}(H^0(E)),\bigwedge^2{\tilde E_2}^*
\boxtimes{\mathcal O}(k-2))$$
$$\rightarrow H^2(\widetilde{S\times S}\times{\mathbb P}(H^0(E)),{\tilde E_2}^*
\boxtimes{\mathcal O}(k-1))$$
is injective.

c) $H^3(\widetilde{S\times S}\times{\mathbb P}(H^0(E)),\bigwedge^3{\tilde E_2}^*
\boxtimes{\mathcal O}(k-3))=0$.

d) The interior product with $\sigma$
$$int(\sigma):H^4(\widetilde{S\times S}\times{\mathbb P}(H^0(E)),\bigwedge^4{\tilde E_2}^*
\boxtimes{\mathcal O}(k-4))$$
$$\rightarrow H^4(\widetilde{S\times S}\times{\mathbb P}(H^0(E)),\bigwedge^3{\tilde E_2}^*
\boxtimes{\mathcal O}(k-3))$$
is injective.}

The conditions a) and c) have been established in proposition \ref{calco}. We now dualize property
b)  as follows : by proposition \ref{calco}, \ref{iti} and \ref{itii},
we have
$$H^2(\widetilde{S\times S}\times{\mathbb P}(H^0(E)),\bigwedge^2{\tilde E_2}^*
\boxtimes{\mathcal O}(k-2))$$
$$=(pr_1^*H^2(S,-L)\oplus pr_2^*H^2(S,-L))\otimes S^{k-2}H^0(S,E)^*,$$

$$H^2(\widetilde{S\times S}\times{\mathbb P}(H^0(E)),{\tilde E_2}^*
\boxtimes{\mathcal O}(k-1))$$
$$=(pr_1^*H^2(S,E^*)\oplus pr_2^*H^2(S,E^*))\otimes S^{k-1}H^0(S,E)^*.$$
Dualizing, we get
$$H^2(\widetilde{S\times S}\times{\mathbb P}(H^0(E)),\bigwedge^2{\tilde E_2}^*
\boxtimes{\mathcal O}(k-2))^*$$
$$=(H^0(S,L)\oplus H^0(S,L))\otimes S^{k-2}H^0(S,E),
$$
$$H^2(\widetilde{S\times S}\times{\mathbb P}(H^0(E)),{\tilde E_2}^*
\boxtimes{\mathcal O}(k-1))^*$$
$$=(H^0(S,E)\oplus H^0(S,E))\otimes S^{k-1}H^0(S,E).
$$
It is then immediate to check that the transpose of the map $int(\sigma)$ is
the map $\wedge\sigma$, so that b) translates into the condition
that
$$\wedge\sigma:(H^0(S,E)\oplus H^0(S,E))\otimes S^{k-1}H^0(S,E)
\rightarrow (H^0(S,L)\oplus H^0(S,L))\otimes S^{k-2}H^0(S,E)$$
is surjective.

Now retracing through the isomorphisms given by proposition \ref{calco},
one checks that the map $\wedge\sigma$ is up to sign equal to the direct sum of two copies of
the composed map
$$\mu:H^0(S,E)\otimes S^{k-1}H^0(S,E)\rightarrow H^0(S,E)\otimes H^0(S,E)\otimes S^{k-2}H^0(S,E)$$
$$
\stackrel{det\otimes id}{\rightarrow }H^0(S,L)\otimes S^{k-2}H^0(S,E).$$

Similarly statement d) dualizes as follows : by proposition \ref{calco}, the 
space $$H^4(\widetilde{S\times S},\bigwedge^4{\tilde E_2}^*\boxtimes{\mathcal 
O}(k-4)) \cong H^4(\widetilde{S\times S},\bigwedge^4{\tilde E_2}^*)\otimes 
S^{k-4}H^0(S,E)^*$$ is dual to $${\rm Ker}\,(H^0(S,L)\otimes H^0(S,L)\rightarrow 
H^0(S,2L))\otimes S^{k-4}H^0(S,E).$$ Next, we know by proposition \ref{calco}, 
\ref{itiv}, that $$H^4(\widetilde{S\times S},\bigwedge^3{\tilde 
E_2}^*\boxtimes{\mathcal O}(k-3)) \cong H^4(\widetilde{S\times 
S},\bigwedge^3{\tilde E_2}^*)\otimes S^{k-3}H^0(S,E)^*$$ admits  a quotient 
which is dual to the direct sum of two copies of $${\rm Ker}\,(H^0(S,E)\otimes 
H^0(S,L)\rightarrow H^0(S,E\otimes L))\otimes S^{k-3}H^0(S,E)).$$ Denoting by 
$Q_{E,L}:={\rm Ker}\,(H^0(S,E)\otimes H^0(S,L)\rightarrow H^0(S,E\otimes L))$, 
$Q_{L,E}:={\rm Ker}\,(H^0(S,L)\otimes H^0(S,E)\rightarrow H^0(S,E\otimes L))$ and 
$Q_{L,L}={\rm Ker}\,(H^0(S,L)\otimes H^0(S,L)\rightarrow H^0(S,2L))$, we have an 
inclusion $$(Q_{L,E}\oplus Q_{E,L}) \otimes S^{k-3}H^0(S,E) \subset 
H^4(\widetilde{S\times S},\bigwedge^3{\tilde E_2}^*\boxtimes{\mathcal 
O}(k-3))^*$$ and to prove d) it suffices to show that the map dual to 
$int(\sigma)$ restricts on this subspace  to a surjection 
$$\wedge\sigma:(Q_{L,E}\oplus Q_{E,L}) \otimes S^{k-3}H^0(S,E) \rightarrow 
Q_{L,L}\otimes S^{k-4}H^0(S,E).$$ But retracing through the isomorphisms of 
proposition \ref{calco} and recalling the definition of $\sigma$, one checks  
easily  that the first component $$\wedge\sigma_1: Q_{L,E}\otimes 
S^{k-3}H^0(S,E) \rightarrow Q_{L,L}\otimes S^{k-4}H^0(S,E)$$ of the map above is 
the following composite 
$$\mu':Q_{L,E}\otimes S^{k-3}H^0(S,E)\subset 
H^0(S,L)\otimes H^0(S,E) \otimes S^{k-3}H^0(S,E) \rightarrow$$ $$H^0(S,L)\otimes 
H^0(S,E)\otimes H^0(S,E)\otimes  S^{k-4}H^0(S,E) $$
$$\stackrel{id\otimes det\otimes 
id}{\rightarrow } H^0(S,L)\otimes H^0(S,L)\otimes S^{k-4}H^0(S,E),$$ which takes 
obviously value in $Q_{L,L}\otimes S^{k-4}H^0(S,E)$, while the second component 
is equal to the first composed with the permutation exchanging factors on both 
sides.

To conclude then that
$$\wedge\sigma:(Q_{L,E}\oplus Q_{E,L}) \otimes S^{k-3}H^0(S,E)
\rightarrow Q_{L,L}\otimes S^{k-4}H^0(S,E)$$
is surjective, it suffices to show that
$$\mu'_-: Q_{L,E}\otimes S^{k-3}H^0(S,E)
\rightarrow Q_{L,L}^-\otimes S^{k-4}H^0(S,E)$$
and 
$$\mu'_+: Q_{L,E}\otimes S^{k-3}H^0(S,E)
\rightarrow Q_{L,L}^+\otimes S^{k-4}H^0(S,E)$$
are surjective, where $Q_{L,L}^+$, (resp. $Q_{L,L}^-$) are the symmetric, resp. antisymmetric
part of $Q_{L,L}$ and $\mu'_+$ (resp. $\mu'_-$) are the composition of
 $\mu'$ with
the projections on the symmetric (resp. antisymmetric) part of $Q_{L,L}$.

In conclusion,  theorem \ref{thmain} will be a consequence of the following
propositions
\begin{prop}\label{prefa} The composed
map
$$\mu:H^0(S,E)\otimes S^{k-1}H^0(S,E)\rightarrow H^0(S,E)\otimes H^0(S,E)
\otimes  S^{k-2}H^0(S,E)$$
$$\stackrel{det}{\rightarrow}H^0(S,L)\otimes S^{k-2}H^0(S,E)$$
is surjective.
\end{prop} 
\begin{prop}\label{prepafa}
a) The  map 
$$\mu'_-:Q_{L,E}\otimes S^{k-3}H^0(S,E)
\rightarrow Q_{L,L}^-\otimes S^{k-4}H^0(S,E)$$ defined above is surjective.

b) The map
$$\mu'_+: Q_{L,E}\otimes S^{k-3}H^0(S,E)
\rightarrow Q_{L,L}^+\otimes S^{k-4}H^0(S,E)$$
defined above is surjective.
\end{prop}
{\bf Proof of proposition \ref{prefa}.} Let $\alpha,\,\beta\in H^0(S,E)$ and
$\gamma\in H^0(S,L)$ such that 
$$\gamma=det(\alpha\wedge\beta).$$
Then we observe first that if $D\subset H^0(S,E)$ is the rank $2$ vector
subspace generated by
$\alpha$ and $\beta$, we have
$$\gamma\otimes S^{k-2}D\subset Im\,\mu$$
since the composite
$$D\otimes S^{k-1}D\rightarrow D\otimes D\otimes S^{k-2}D\rightarrow
 \bigwedge^2D
\otimes S^{k-2}D$$
is surjective.

Recall now that the map $det$ determines a morphism
$$d:G_2\rightarrow {\mathbb P}H^0(S,L)$$
which is surjective and finite since both spaces are of the same dimension
$2k$. The fiber $d^{-1}(\gamma)$ is then a finite subscheme
$Z_\gamma\subset G_2$ which is the complete intersection of
a space $W$ of 
hyperplane sections of the Grassmannian $G_2$.

Now by the above observation, and since $d$ is surjective, 
 it suffices to show that
the subspaces $S^{k-2}D$ for $D\in Z_\gamma$
generate $S^{k-2}H^0(S,E)$. If we dualize, this is equivalent
 to say that the dual
map
$$S^{k-2}H^0(S,E)^*\rightarrow \oplus_{D\in Z_\gamma}S^{k-2}D^*$$
is injective. But this map identifies to the restriction
$$H^0(G_2,S^{k-2}{\mathcal E})\rightarrow 
H^0(Z_\gamma,S^{k-2}{\mathcal E}_{\mid Z_\gamma}),$$
at least for a reduced $Z_\gamma$, which will be 
the case for a generic $\gamma$.

Hence it suffices to show that
\begin{eqnarray}
\label{encorevan}
H^0(G_2,S^{k-2}{\mathcal E}\otimes{\mathcal I}_{Z_\gamma})=0.\end{eqnarray}
Now we use the Koszul resolution
$$0\rightarrow \bigwedge^{2k}W\otimes{\mathcal L}^{-2k}\rightarrow \ldots
\rightarrow W\otimes{\mathcal L}^{-1}\rightarrow {\mathcal I}_{Z_\gamma}\rightarrow 0$$

The vanishing (\ref{encorevan}) will then follow from the vanishing
$$H^i(G_2,S^{k-2}{\mathcal E}\otimes{\mathcal L}^{-i-1}),\,i=0,\,2k-1$$
which is proved in proposition \ref{appen} of the appendix.
\cqfd

{\bf Proof of proposition \ref{prepafa}, a).}  Notice first that
the natural composed map
$$\bigwedge^3H^0(S,E)\rightarrow \bigwedge^2H^0(S,E)\otimes H^0(S,E)$$
$$\stackrel{det\otimes id}{\rightarrow } H^0(S,L)\otimes H^0(S,E)$$
has its image contained in $Q_{L,E}$.
Hence it suffices to show that the following composite
$$\mu'':\bigwedge^3H^0(S,E)\otimes S^{k-3}H^0(S,E)\rightarrow \bigwedge^2H^0(S,E)\otimes H^0(S,E)
\otimes S^{k-3}H^0(S,E)$$
$$
\stackrel{det\otimes\mu}{\rightarrow }H^0(S,L)\otimes H^0(S,L)\otimes S^{k-4}H^0(S,E)
\rightarrow \bigwedge^2H^0(S,L)\otimes S^{k-4}H^0(S,E)$$
is surjective.

Now note that for $\alpha_1,\,\alpha_2,\,\alpha_3\in H^0(S,E)$
\begin{eqnarray}
\label{eqdim}\mu''(\alpha_1\wedge\alpha_2\wedge\alpha_3\otimes \alpha_3^{k-3})=
2(k-3)det(\alpha_2\wedge\alpha_3)\wedge det(\alpha_1\wedge\alpha_3)\otimes\alpha_3^{k-4}.
\end{eqnarray}
Fix now $\gamma\in H^0(S,L)$ and consider the set of couples
$(\alpha_1,\alpha_3)$ such that 
$$det(\alpha_1\wedge\alpha_3)=\gamma.$$
 For any
$\alpha_2$ and any 
such $(\alpha_1,\alpha_3)$, we have
$$\mu''(\alpha_1\wedge\alpha_2\wedge\alpha_3\otimes \alpha_3^{k-3})=
2(k-3)det(\alpha_2\wedge\alpha_3)\wedge \gamma\otimes\alpha_3^{k-4}.$$
Note that the vector $\alpha_3$ for such pairs takes arbitrary value in some
of the lines $D\in Z_\gamma$, where the notations are as in the previous proposition.

Now we have the map
$$\mu''':H^0(S,E)\otimes S^{k-3}H^0(S,E)\rightarrow H^0(S,L)\otimes S^{k-4}H^0(S,E)$$
analogous to $\mu$ and the formula above
shows that
$$\mu''(\alpha_1\wedge\alpha_2\wedge\alpha_3\otimes
 \alpha_3^{k-3})=2\gamma\wedge\mu'''(\alpha_2\otimes\alpha_3^{k-3}).$$
 With the same proof as in the previous proposition, one shows now that
 the $S^{k-3}D,\,D\in Z_\gamma$ generate $S^{k-3}H^0(S,E)$ and that $\mu'''$ is surjective. Hence
 the $\alpha_2\otimes\alpha_3^{k-3},\,\alpha_3\in D,\,D\in Z_\gamma$
 generate $H^0(S,E)\otimes S^{k-3}H^0(S,E)$ and
 the $\mu'''(\alpha_2\otimes\alpha_3^{k-3}),\,
 \alpha_3\in D,\,D\in Z_\gamma$ generate by the surjectivity of $\mu'''$ the space
 $H^0(S,L)\otimes S^{k-4}H^0(S,E)$. Hence
 $Im\,\mu''$ contains $\gamma\wedge H^0(S,L)\otimes S^{k-4}H^0(S,E)$, and since $\gamma$
 was generic, we conclude that $\mu''$ is surjective.

\cqfd
{\bf Proof of proposition \ref{prepafa}, b).}
We want to prove that
$$\mu'_+:Q_{L,E}\otimes S^{k-3}H^0(S,E)\rightarrow 
Q_{L,L}^+\otimes S^{k-4}H^0(S,E)$$
is surjective. Denote similarly, for $C$ a generic member of $\mid L\mid$,
$$Q_{K_C,E}:={\rm Ker}\,(H^0(C,K_C)\otimes H^0(C,E_{\mid C})\rightarrow 
H^0(C,E\otimes K_C)),$$
$$Q_{K_C,K_C}^+:={\rm Ker}\,(S^2H^0(C,K_C)\rightarrow 
H^0(C, K_C^{\otimes2})).$$
Then we can define 
$$\mu'_{+,C}:Q_{K_C,E}\otimes S^{k-3}H^0(C,E_{\mid C})
\rightarrow Q_{K_C,K_C}^+\otimes S^{k-4}H^0(C,E_{\mid C})$$
as the composite
$$Q_{K_C,E}\otimes S^{k-3}H^0(C,E_{\mid C})\subset H^0(C,K_C)\otimes H^0(C,E_{\mid C})
\otimes S^{k-3}H^0(C,E_{\mid C})$$
$$\rightarrow H^0(C,K_C)\otimes H^0(C,E_{\mid C})\otimes H^0(C,E_{\mid C})\otimes S^{k-4}H^0(C,E_{\mid C})
\stackrel{id\otimes det\otimes id}{\rightarrow}$$
 $$H^0(C,K_C)\otimes H^0(C,K_C)\otimes S^{k-4}H^0(C,E_{\mid C})
\rightarrow S^2H^0(C,K_C)\otimes S^{k-4}H^0(C,E_{\mid C}).$$
Now the restriction map
$H^0(S,E)\rightarrow H^0(C,E_{\mid C})$ is an isomorphism, and
the restriction map $H^0(S,L)\rightarrow H^0(C,K_C)$
is surjective with kernel $\sigma$, the defining equation of
$C$. Hence
the restrictions induce a surjection
$$Q_{L,E}\rightarrow  Q_{K_C,E}$$
and an isomorphism
$$Q_{L,L}^+\cong Q_{K_C,K_C}^+,$$
and it suffices to show that
$\mu'_{+,C}$ is surjective. A fortiori it suffices to show that the composite
$$\mu'_{C}:Q_{K_C,E}\otimes S^{k-3}H^0(C,E_{\mid C})
\subset H^0(C,K_C)\otimes H^0(C,E_{\mid C})
\otimes S^{k-3}H^0(C,E_{\mid C})$$
$$\rightarrow H^0(C,K_C)\otimes H^0(C,E_{\mid C})
\otimes H^0(C,E_{\mid C})\otimes S^{k-4}H^0(C,E_{\mid C})$$
$$
\stackrel{id\otimes det\otimes id}{\rightarrow}
 H^0(C,K_C)\otimes H^0(C,K_C)\otimes S^{k-4}H^0(C,E_{\mid C})$$
 which takes value in $Q_{K_C,K_C}:={\rm Ker}\,(H^0(C,K_C)^{\otimes2}\rightarrow 
H^0(C, K_C^{\otimes2}))$, is surjective on this last space.

Let us now consider the following diagram of exact sequences
$$\begin{matrix}
&0&\rightarrow &Q_{K_C,E}\otimes S^{k-3}H^0(C,E_{\mid C})&\rightarrow &
H^0(C,K_C)\otimes H^0(C,E_{\mid C})
\otimes S^{k-3}H^0(C,E_{\mid C})\\
&&&\mu'_C\downarrow&& id\otimes\mu_C\downarrow&\\
&0&\rightarrow &Q_{K_C,K_C}\otimes S^{k-4}H^0(C,E_{\mid C})&\rightarrow &
H^0(C,K_C)\otimes H^0(C,K_C)\otimes S^{k-4}H^0(C,E_{\mid C})
\end{matrix}
$$
$$\begin{matrix}
&\rightarrow &H^0(C,E\otimes K_C)\otimes S^{k-3}H^0(C,E_{\mid C})&
\rightarrow &0\\
&&\mu_{C,K_C}\downarrow&&\\
&\rightarrow &
H^0(C, K_C^{\otimes2})\otimes S^{k-4}H^0(C,E_{\mid C})&\rightarrow &0
\end{matrix},
$$

 where
 the vertical maps
$\mu_C$ and $\mu_{C,K_C}$ are defined in a way similar
to $\mu$ e.g
$\mu_C$ is the composite
$$H^0(C,E_{\mid C})
\otimes S^{k-3}H^0(C,E_{\mid C})\subset H^0(C,E_{\mid C})\otimes H^0(C,E_{\mid C})
\otimes S^{k-4}H^0(C,E_{\mid C})$$
$$\stackrel{det\otimes id}{\rightarrow }
H^0(C,K_C)\otimes S^{k-4}H^0(C,E_{\mid C}),$$
and $\mu_{C,K_C}$ is defined similarly with a twist by $K_C$.

One checks easily the surjectivity of the multiplication maps on the left.

The proof of  proposition \ref{prefa} shows as well that 
$\mu_C$ is surjective, as is $\mu_{C,K_C}$ by the commutativity of the diagram above.
Hence the surjectivity of $\mu'_C$ will follow by diagram chasing from the surjectivity
of the induced multiplication map
\begin{eqnarray}
\label{mullun}
H^0(C,K_C)\otimes {\rm Ker}\,\mu_C\rightarrow {\rm Ker}\,\mu_{C,K_C}.
\end{eqnarray}
In what follows we will use again the notation $H^0(E)$ for $H^0(S,E)=H^0(C,E_{\mid C})$.
Define the vector bundle ${\mathcal Q}$ on $C$ as the kernel of the surjective
composite  morphism of vector bundles
$$S^{k-3}H^0(E)\otimes E\subset S^{k-4}H^0(E)\otimes
H^0(E)\otimes E\stackrel{id\otimes det}{\rightarrow }S^{k-4}H^0(E)\otimes K_C.$$
Then we clearly have
$${\rm Ker}\,\mu_C=H^0(C,{\mathcal Q}),\, {\rm Ker}\,\mu_{C,K_C}=H^0(C,{\mathcal Q}\otimes K_C)$$
so that the surjectivity of the map (\ref{mullun})
is equivalent to the surjectivity of the multiplication map
\begin{eqnarray}
\label{mullunbun}H^0(C,{\mathcal Q})\otimes H^0(C,K_C)\rightarrow H^0(C,{\mathcal Q}\otimes K_C).
\end{eqnarray}
Now we proceed as follows : let $\sigma\in H^0(S,L)$ be the defining equation for
$C$. Recall the finite reduced subscheme $Z_\sigma=d^{-1}(\sigma)\subset G_2$
made of the rank $2$ vector subspaces $D$ of $H^0(S,E)$ such that $det\,D=\sigma$.
For each such $D$, there is a line subbundle $L_D$ of $E$ on $C$,
of degree $k+1$ with two sections without common zeroes (see section \ref{sec1}). The
space $D$ identifies naturally to $H^0(C,L_D)$.

Clearly the image of the inclusion
$$S^{k-3}H^0(C,L_D)\otimes L_D\subset S^{k-3}H^0(E)\otimes E$$
is contained in ${\mathcal Q}$.

Let now
$${\mathcal N}:=\oplus_{D\in Z_\sigma}S^{k-3}D\otimes L_D.$$
Then by the observation above we have a morphism
$$\alpha: {\mathcal N}\rightarrow {\mathcal Q}.$$
The   surjectivity of
\ref{mullunbun} will follow from the following three lemmas :
\begin{lemm}\label{lesurlun} The morphism $\alpha$ is surjective.
\end{lemm}
Denoting ${\mathcal M}:={\rm Ker}\,\alpha$ we also prove
\begin{lemm}\label{secsurlun} The vector bundle ${\mathcal M}$ is generated by its sections.
\end{lemm}
\begin{lemm} \label{sur-x}The space $H^0(C,{\mathcal M})$ is generated by
the subspaces $H^0(C,{\mathcal M}(-x)),\,x\in C$.
\end{lemm}
We explain first how these three lemmas imply our result. Using the exact sequence
$$0\rightarrow {\mathcal M}\rightarrow {\mathcal N}\rightarrow {\mathcal Q}
\rightarrow 0$$
given by lemma \ref{lesurlun}, we see that the map (\ref{mullunbun})
will be surjective if 
the multiplication map
$$H^0(C,{\mathcal N})\otimes H^0(C,K_C)\rightarrow H^0(C,{\mathcal N}\otimes K_C)$$
is surjective, and $H^1(C,{\mathcal M}\otimes K_C)=0$.

The first condition is easy to check. Indeed ${\mathcal N}$ is a direct sum of line bundles
$L_D$ corresponding to $g_{k+1}^1$'s on $C$, and the result is easy to prove for each of them.
As for the second condition, it is equivalent
to $H^0(C,{\mathcal M}^*)=0$ by Serre's duality. But since ${\mathcal M}$
is generated by sections by lemma \ref{secsurlun}, we have an inclusion
$$H^0(C,{\mathcal M}^*)\subset H^0(C,{\mathcal M})^*.$$
The image of this inclusion obviously vanishes on each 
subspace $H^0(C,{\mathcal M}(-x))$, hence it must be $0$ since we know by lemma
\ref{sur-x} that these subspaces generate $H^0(C,{\mathcal M})$.
\cqfd
To conclude the proof of   \ref{prepafa},b) it remains only to prove
these three lemmas.

\vspace{0,5cm}

{\bf Proof of lemma \ref{lesurlun}.} First of all we note that the bundle ${\mathcal Q}$
is generated by its sections, since there is a natural surjection
$$S^{k-2}H^0(E)\otimes{\mathcal O}_C\rightarrow {\mathcal Q}\rightarrow 0.$$
Hence it suffices to show that the map
$$H^0(C,{\mathcal N})\rightarrow H^0(C,{\mathcal Q})$$
is surjective.

But by definition 
$$H^0(C,{\mathcal Q})={\rm Ker}\,(H^0(E)\otimes S^{k-3}H^0(E)\stackrel{\mu_C}
{ \rightarrow}
 H^0(C,K_C)\otimes S^{k-4}H^0(E)) $$
 and $$H^0(C,{\mathcal N})=\oplus_{D\in Z_\sigma}D\otimes S^{k-3}D.$$
 Hence we need to show that the sequence
 $$\oplus_{D\in Z_\sigma}D\otimes S^{k-3}D\rightarrow 
 H^0(E)\otimes S^{k-3}H^0(E)\stackrel{\mu_C}{ \rightarrow}
 H^0(C,K_C)\otimes S^{k-4}H^0(E)$$
 is exact at the middle. Again this will
  follow from a cohomological computation on the 
 Grassmannian $G_2$. Indeed, the notations being as in the
 proof of Propositions \ref{newform} and  \ref{prefa},
  the sequence above dualizes as
 \begin{eqnarray}
 \label{g2ex}
  I_{Z_\sigma}({\mathcal L})\otimes S^{k-4}H^0(G_2,{\mathcal E})
 \rightarrow H^0(G_2,{\mathcal E})\otimes S^{k-3}H^0(G_2,{\mathcal E})
 \nonumber\\
 \rightarrow H^0({\mathcal E}\otimes S^{k-3}{\mathcal E}_{\mid Z_\sigma}),
 \end{eqnarray}
  where 
 the map
 $$I_{Z_\sigma}({\mathcal L})\otimes S^{k-4}H^0(G_2,{\mathcal E})
 \rightarrow H^0(G_2,{\mathcal E})\otimes S^{k-3}H^0(G_2,{\mathcal E})$$
 is composed of the inclusion
 $$I_{Z_\sigma}({\mathcal L})\otimes S^{k-4}H^0(G_2,{\mathcal E})
 \subset H^0(G_2,{\mathcal L})\otimes S^{k-4}H^0(G_2,{\mathcal E})$$
 $$
 \cong \wedge^2H^0(G_2,{\mathcal E})\otimes S^{k-4}H^0(G_2,{\mathcal E})$$
 and of the (Koszul) map
 $$\wedge^2H^0(G_2,{\mathcal E})\otimes S^{k-4}H^0(G_2,{\mathcal E})
 \rightarrow H^0(G_2,{\mathcal E})\otimes S^{k-3}H^0(G_2,{\mathcal E}).$$
 One checks easily that
 $H^0(G_2,{\mathcal E})\otimes S^{k-3}H^0(G_2,{\mathcal E})
 \cong H^0(G_2,{\mathcal E}\otimes S^{k-3}{\mathcal E})$. Hence the kernel in the middle identifies to
 $H^0(G_2,{\mathcal E}\otimes S^{k-3}{\mathcal E}\otimes{\mathcal I}_{Z_\sigma})$.
 Furthermore $S^{k-4}H^0(G_2,{\mathcal E})\cong H^0(G_2,S^{k-4}{\mathcal E})$
 identifies to $H^0(G_2,{\mathcal E}\otimes S^{k-3}{\mathcal E}\otimes{\mathcal L}^{-1})$
 via the (Koszul) inclusion
 $$S^{k-4}{\mathcal E}\otimes{\mathcal L}=S^{k-4}{\mathcal E}\otimes\bigwedge^2{\mathcal E}
 \subset {\mathcal E}\otimes S^{k-3}{\mathcal E}.$$
 
 Hence the exactness at the middle of the sequence \ref{g2ex}
 will follow from 
 the surjectivity of the multiplication map
 \begin{eqnarray}
 \label{plouf}
 H^0(G_2,{\mathcal E}\otimes S^{k-3}{\mathcal E}
 \otimes{\mathcal L}^{-1})\otimes I_{Z_\sigma}({\mathcal L})
 \rightarrow H^0(G_2,{\mathcal E}\otimes S^{k-3}{\mathcal E}\otimes{\mathcal I}_{Z_\sigma})
 .
 \end{eqnarray}
 
 Now let $W:=I_{Z_\sigma}({\mathcal L})$. The Koszul resolution of ${\mathcal I}_{Z_\sigma}$
 $$0\rightarrow \bigwedge^{2k}W\otimes{\mathcal L}^{-2k}\rightarrow \ldots
 \rightarrow W\otimes{\mathcal L}^{-1}\rightarrow {\mathcal I}_{Z_\sigma}\rightarrow 0$$
 twisted by ${\mathcal E}\otimes S^{k-3}{\mathcal E}$ shows that the surjectivity
 of
 (\ref{plouf})
 will hold if we know that
 $$H^{i}(G_2,{\mathcal E}\otimes S^{k-3}{\mathcal E}\otimes{\mathcal L}^{-i-1})=0,\,1\leq i<2k.$$
 Since we have the exact sequence
 $$0\rightarrow S^{k-4}{\mathcal E}\otimes{\mathcal L}\rightarrow {\mathcal E}\otimes S^{k-3}{\mathcal E}
 \rightarrow S^{k-2}{\mathcal E}\rightarrow 0,$$
 it suffices  to know that
 $$H^{i}(G_2, S^{k-4}{\mathcal E}\otimes{\mathcal L}^{-i})=0,\,1\leq i<2k,$$
 and
$$ H^{i}(G_2, S^{k-2}{\mathcal E}\otimes{\mathcal L}^{-i-1})=0,\,1\leq i<2k.$$
 This is proved in Proposition \ref{appen}.
 \cqfd
{\bf Proof of lemma \ref{secsurlun}.}
The bundles ${\mathcal N}$ and ${\mathcal Q}$ are generated by global sections.
To prove that ${\mathcal M}$ is generated by global sections, it suffices to prove that for
any
$x\in C$, the restriction map $H^0(C,{\mathcal N}(-x))\rightarrow H^0(C,{\mathcal Q}(-x))$
is surjective. For each $g_{k+1}^1\,\,L_D$ on $C$, denote
by $\sigma_{D,x}\in H^0(C,L_D)\cong D$ a generator for
$H^0(C,L_D(-x))$. We need to show the exactness
of the  sequence
\begin{eqnarray}
\label{encexlun}\oplus_{D\in Z_{\sigma}}\sigma_{D,x}\otimes S^{k-3}D
\rightarrow H^0(C,E(-x))\otimes S^{k-3}H^0(E)\nonumber \\
\stackrel{\mu_C}{\rightarrow} 
H^0(C,K_C(-x))\otimes S^{k-4}H^0(E).
\end{eqnarray}
(Indeed, by definition,  $\oplus_{D\in Z_\sigma}\sigma_{D,x}\otimes S^{k-3}D$
identifies to $H^0(C,{\mathcal N}(-x))$ and
$$ {\rm Ker}\,(H^0(C,E(-x))\otimes S^{k-3}H^0(E)\stackrel{\mu_C}{\rightarrow}
 H^0(C,K_C(-x))\otimes S^{k-4}H^0(E))$$
identifies to $H^0(C,{\mathcal Q}(-x))$.)

Denote by $K_x\subset H^0(E)$ the subspace $H^0(C,E(-x))$. Note that via the identification
$H^0(C,L_D)=D\subset H^0(E)$,
 $\sigma_{D,x}$ becomes a generator of the $1$-dimensional vector space
$D\cap K_x$. Furthermore, $K_x$ being of codimension
$2$ in $H^0(E)$ determines a section
$\tau_x\in\bigwedge^2H^0(E)^*$ up to a coefficient. Clearly 
$\tau_x$ belongs to $ H^0(C,K_C)^*\subset \bigwedge^2H^0(E)^*$ 
and identifies also  to the
linear form on $H^0(C,K_C)$ defining $H^0(C,K_C(-x))$. Let
$G_x\subset G_2$ be the hyperplane section defined by $\tau_x$. The scheme $Z_\sigma$ is a complete
intersection of hyperplane sections of 
$G_x$. The variety $G_x$ admits a desingularization
$P_x\stackrel{g}{\rightarrow }G_x$ defined as
$$P_x=\{(u,\Delta)\in{\mathbb P}(K_x)\times G_2,\,u\in\Delta\cap K_x\}.$$
Note that if
$$\begin{matrix} &P&\stackrel{g}{\rightarrow }&G_2&\\
&f\downarrow&&&\\
&{\mathbb P}(H^0(E))&&&
\end{matrix}
$$
is the incidence variety, $P_x$ can also be defined as $f^{-1}({\mathbb P}(K_x))\subset P$.

Each line $D$ parametrized by $Z_\sigma$ meets 
$K_x$ along a one dimensional vector space, because $L_D$
has no base-point, so that
$H^0(L_D(-x))\cong D\cap K_x$ is $1$-dimensional.
 It follows that the scheme $Z_\sigma$ can also be seen
as the complete intersection in $P_x$ of hypersurfaces in $\mid g^*{\mathcal L}\mid$.

We now dualize the sequence (\ref{encexlun}). The space $H^0(C,K_C(-x))$
admits for dual the space
$W\subset H^0(P_x,g^*{\mathcal L})$ defining $Z_\sigma\subset P_x$. The vector space 
$<\sigma_{D,x}>^*$ identifies clearly to the fiber of the line bundle
$H_x:=f^*{\mathcal O}_{{\mathbb P}(K_x)}(1)$ at the point $D\in Z_\sigma$. 
Hence our sequence dualizes as
\begin{eqnarray}
\label{dualenclun}
W\otimes S^{k-4}H^0(E)^*
\rightarrow f^*H^0({\mathbb P}(K_x),{\mathcal O}(1))
\otimes H^0(P_x,g^*S^{k-3}{\mathcal E}) \nonumber \\
\rightarrow H^0(S^{k-3}{\mathcal E}\otimes H_x{\mid Z_\sigma}).
\end{eqnarray}
The second space in this sequence  is easily shown to identify to
$H^0(P_x,g^*S^{k-3}{\mathcal E}\otimes H_x)$, so that the kernel at the middle
is equal to
$$H^0(P_x,g^*S^{k-3}{\mathcal E}\otimes H_x\otimes{\mathcal I}_{Z_\sigma}).$$
The first map in (\ref{dualenclun}) is induced by the isomorphism
$$S^{k-4}H^0(E)^*\cong H^0(P_x,g^*S^{k-4}{\mathcal E}),$$
 the multiplication
$$W\otimes H^0(P_x,g^*S^{k-4}{\mathcal E})\rightarrow 
H^0(P_x,g^*(S^{k-4}{\mathcal E}\otimes{\mathcal L})\otimes
{\mathcal I}_{Z_\sigma})$$
and 
the composed bundle map
$$g^*S^{k-4}{\mathcal E}\otimes{\mathcal L}\rightarrow g^*S^{k-3}{\mathcal E}\otimes{\mathcal E}
\rightarrow g^*S^{k-3}{\mathcal E}\otimes H_x,$$
where the last map is induced by the  natural surjective
map
$g^*{\mathcal E}\rightarrow H_{x}$. 
 
The exactness of (\ref{dualenclun}) will then follow from the surjectivity of
\begin{eqnarray}
\label{enkolun}
W\otimes H^0(P_x,g^*(S^{k-3}{\mathcal E}\otimes{\mathcal L}^{-1})\otimes H_x)
\rightarrow H^0(P_x,g^*S^{k-3}{\mathcal E}\otimes H_x\otimes{\mathcal I}_{Z_\sigma})
\end{eqnarray}
and from the equality
\begin{eqnarray}
\label{eqmi}
H^0(P_x,g^*(S^{k-3}{\mathcal E}\otimes{\mathcal L}^{-1})\otimes H_x)=
H^0(P_x,g^*S^{k-4}{\mathcal E}).
\end{eqnarray}
This last equality is proved as follows : on $P_x$ we have the exact sequence
$$0\rightarrow g^*{\mathcal L}\otimes H_x^{-1}\rightarrow g^*{\mathcal E}\rightarrow H_x\rightarrow 0,$$
which gives
$$0\rightarrow g^*S^{k-4}{\mathcal E}\otimes g^*{\mathcal L}\otimes H_x^{-1}
\rightarrow g^*S^{k-3}{\mathcal E}\rightarrow H_x^{k-3}\rightarrow 0.$$
Tensoring this with $H_x\otimes{\mathcal L}^{-1}$ we get
$$0\rightarrow g^*S^{k-4}{\mathcal E}\rightarrow g^*(S^{k-3}{\mathcal E}
\otimes{\mathcal L}^{-1})\otimes H_x
\rightarrow g^*{\mathcal L}^{-1}\otimes
H_x^{k-2} \rightarrow 0.$$
But the right hand side has no non-zero sections since it is of negative degree
on the fibers of $f$.  Hence the equality
(\ref{eqmi}). 

Since $Z_\sigma\subset P_x$ is the complete intersection of the space $W$ of sections of
$g^*{\mathcal L}$, we have a Koszul resolution of ${\mathcal I}_{Z_\sigma}$, which takes the form
$$0\rightarrow \bigwedge^{2k-1}W\otimes g^*{\mathcal L}^{-2k-1}\rightarrow 
\ldots\rightarrow  W\otimes g^*{\mathcal L}^{-1}\rightarrow{\mathcal I}_{Z_\sigma}\rightarrow 0.$$
We can tensor it with $g^*S^{k-3}{\mathcal E}\otimes H_x$, and the surjectivity of the
map (\ref{enkolun}) will follow from the following vanishing
\begin{eqnarray}
\label{vanvanlun} H^i(P_x,g^*(S^{k-3}{\mathcal E}\otimes {\mathcal L}^{-i-1})\otimes
 H_x)=0,\,1\leq i<2k-1=dim\,P_x.
\end{eqnarray}
 Recall now that $P_x\subset P$ is the complete intersection of two sections of
 $H=f^*{\mathcal O}_{{\mathbb P}(H^0(E))}(1)$, with $H_x=H_{\mid P_x}$.
 The vanishing (\ref{vanvanlun}) will then  follow
 from
 $$H^i(P,g^*(S^{k-3}{\mathcal E}\otimes{\mathcal L}^{-i-1})\otimes
 H)=0,\,1\leq i<2k-1$$
 $$ H^{i+1}(P,g^*(S^{k-3}{\mathcal E}
 \otimes {\mathcal L}^{-i-1}))=0,\,1\leq i<2k-1$$
 $$H^{i+2}(P,g^*(S^{k-3}{\mathcal E}\otimes
 {\mathcal L}^{-i-1})\otimes
 H^{-1})=0,\,1\leq i<2k-1.$$
 The second equality follows immediately from the proposition
 \ref{appen}, and the third is obvious since $H^{-1}$ has no cohomology on the fibers of
 $g$. The first equality is proven as follows :
 we have
 $$H^i(P,g^*S^{k-3}{\mathcal E}\otimes{\mathcal L}^{-i-1})
\otimes H )=H^i(G_2,S^{k-3}{\mathcal E}\otimes{\mathcal E}
 \otimes{\mathcal L}^{-i-1}),$$
 since $R^0g_*H={\mathcal E}$. Now we have the exact sequence
 on $G_2$
 $$0\rightarrow S^{k-4}{\mathcal E}\otimes{\mathcal L}
 \rightarrow S^{k-3}{\mathcal E}\otimes{\mathcal E}\rightarrow S^{k-2}{\mathcal E}\rightarrow 0.$$
 Hence the needed equality will follow from
 the vanishings
 $$H^i(G_2,S^{k-4}{\mathcal E}
 \otimes{\mathcal L}^{-i})=0,$$
 $$H^i(G_2,S^{k-2}{\mathcal E}
 \otimes{\mathcal L}^{-i-1})=0,$$
 for $1\leq i<2k-1$, which are proved in proposition (\ref{appen}).
 Hence lemma  \ref{secsurlun} is proven.
 
\cqfd
{\bf Proof of lemma \ref{sur-x}.}
Let $x_1,\ldots,x_{2k-1}$ be points of $C$ in general position. We will show that
the natural map
\begin{eqnarray}
\label{dersurlun}
\oplus_iH^0(C,{\mathcal M}(-x_i))\rightarrow H^0(C,{\mathcal M})
\end{eqnarray}
is surjective.

Recall that 
$$H^0(C,{\mathcal M})={\rm Ker}\,\oplus_{D\in Z_\sigma}S^{k-3}D\otimes D\rightarrow S^{k-3}H^0(E)\otimes H^0(E).$$
It follows from this, using the identifications
$$H^0(E)^*=H^0(G_2,{\mathcal E}),\,D^*\cong{\mathcal E}_D$$
 that
$$H^0(C,{\mathcal M})^*=Coker\,H^0(G_2,S^{k-3}{\mathcal E}\otimes {\mathcal E})
\rightarrow H^0(S^{k-3}{\mathcal E}\otimes {\mathcal E}_{\mid Z_\sigma})$$
$$=
H^1(G_2,S^{k-3}{\mathcal E}\otimes {\mathcal E}\otimes{\mathcal I}_{Z_\sigma}).$$
Similarly
$$H^0(C,{\mathcal M}(-x_i))={\rm Ker}\,\oplus_{D\in Z_\sigma} S^{k-3}D\otimes\sigma_{D,x_i}
\rightarrow S^{k-3}H^0(E)\otimes K_{x_i}$$
which, with the notations of the previous proof, dualizes to
$$H^0(C,{\mathcal M}(-x_i))^*={\rm Coker}\,(H^0(P_{x_i},g^*S^{k-3}{\mathcal E}\otimes H_{x_i})
\rightarrow H^0(Z_\sigma,S^{k-3}{\mathcal E}\otimes H_{x_i}))$$
$$=
H^1(P_{x_i},g^*S^{k-3}{\mathcal E}\otimes H_{x_i}\otimes{\mathcal I}_{Z_\sigma}),
$$
where we view $Z_\sigma$ as a subscheme of $P_{x_i}$ as well.
Hence we have to show that the natural  map (induced by the morphism
$g^*{\mathcal E}\rightarrow H_{x_i}$ on $P_{x_i}$)
\begin{eqnarray}
\label{sannom}
H^1(G_2,S^{k-3}{\mathcal E}\otimes{\mathcal E}\otimes{\mathcal I}_{Z_\sigma})\rightarrow 
\oplus_iH^1(P_{x_i},g^*S^{k-3}{\mathcal E}\otimes H_{x_i}\otimes{\mathcal I}_{Z_\sigma})
\end{eqnarray}
is injective.

Let $R\subset G_2$ be the curve which is the complete intersection of the
sections $\sigma_{x_i}\in H^0(G_2,{\mathcal L})$. We have first

{\bf Fact.}{\it  The restriction map
$$H^0(G_2,S^{k-3}{\mathcal E}\otimes{\mathcal E})\rightarrow 
H^0(R,S^{k-3}{\mathcal E}\otimes{\mathcal E}_{\mid R})$$
is surjective.}

Using the Koszul resolution of ${\mathcal I}_D$ this is obtained by application of the proposition
\ref{appen}.
\cqfd

>From this we conclude that the restriction
map 
$$H^1(G_2,S^{k-3}{\mathcal E}\otimes{\mathcal E}\otimes{\mathcal I}_{Z_\sigma})
\rightarrow H^1(R,S^{k-3}{\mathcal E}\otimes{\mathcal E}_{\mid R}\otimes{\mathcal I}_{Z_\sigma})$$
is injective.

Consider now the fibered product
$$\tilde{R}:=P_{x_1}\times_{G_2}\times\ldots\times_{G_2}P_{x_{2k-1}}.$$
Denote by $\tilde g:\tilde R\rightarrow R\subset G_2$ the natural morphism.
One shows easily that the curve $\tilde R$ is isomorphic to $R$ excepted over
the intersection of $R$ with the Grassmannian of lines in ${\mathbb P}(K_{x_i})$
for some $i$. Here $R$ has nodes, which are replaced in $\tilde{R}$ by lines.
 
(This fact is obviously true set theoretically, and is proved scheme theoretically by the
computation of the canonical bundles, which gives :
$$K_{\tilde R}=\tilde g^* K_R.)$$
The zero set
$Z_\sigma$ is supported away of this singular locus.
For each $i$ we have a natural restriction map
$$H^1(P_{x_i},g^*S^{k-3}{\mathcal E}\otimes H_{x_i}\otimes{\mathcal I}_{Z_\sigma})
\rightarrow H^1(\tilde R,\tilde g^*S^{k-3}{\mathcal E}\otimes H_{x_i}\otimes{\mathcal I}_{Z_\sigma}),$$
since 
$\tilde R=P_{x_1}\times_{G_2}\ldots\times_{G_2} P_{x_{2k-1}}$
admits a natural morphism to $P_{x_i}$.
Next we have by the above description of
$\tilde R$ an isomorphism
$$H^1(R,S^{k-3}{\mathcal E}\otimes{\mathcal E}_{\mid R}\otimes{\mathcal I}_{Z_\sigma})
\cong H^1(\tilde R,\tilde g^*S^{k-3}{\mathcal E}\otimes{\mathcal E}\otimes{\mathcal I}_{Z_\sigma})$$
and it follows that the injectivity of the map
(\ref{sannom}) will be a consequence of the injectivity of the 
 map
\begin{eqnarray}
\label{sannomnom}
H^1(\tilde R,\tilde g^*S^{k-3}{\mathcal E}\otimes{\mathcal E}\otimes{\mathcal I}_{Z_\sigma})
\rightarrow \oplus_iH^1(\tilde R,\tilde g^*S^{k-3}{\mathcal E}
\otimes H_{x_i}\otimes{\mathcal I}_{Z_\sigma})
\end{eqnarray}
induced by the morphisms
$\tilde g^*{\mathcal E}\rightarrow H_{x_i}$ on $\tilde R$.
Recall now that $Z_\sigma\subset R$ is defined by
a section of ${\mathcal L}$ so that similarly
$Z_\sigma\subset \tilde R$ is defined by
a section of $\tilde g^*{\mathcal L}$. Hence we have
$${\mathcal I}_{Z_\sigma}\cong\tilde g^*{\mathcal L}^{-1}.$$
Furthermore
$$K_{\tilde R}=\tilde g^* K_R=\tilde g^*({K_{G_2}}_{\mid R}\otimes{\mathcal L}^{2k-1})=
\tilde g^*{\mathcal L}^{k-3}.$$
Hence the map (\ref{sannomnom})
 dualizes by Serre's duality as the map
\begin{eqnarray}
\label{saumon}
\oplus_iH^0(\tilde R,\tilde g^*S^{k-3}{\mathcal E}^*\otimes H_{x_i}^*\otimes
\tilde g^*{\mathcal L}\otimes \tilde g^*{\mathcal L}^{k-3})
\\
\rightarrow H^0(\tilde R,\tilde g^*(S^{k-3}{\mathcal E}^*\otimes {\mathcal E}^*)\otimes
\tilde g^*{\mathcal L}
\otimes \tilde g^*{\mathcal L}^{k-3})
 \end{eqnarray}
 given by the inclusions
 $H_{x_i}^*\subset\tilde g^*{\mathcal E}^*$ on $\tilde R$. Since
 $det\,{\mathcal E}={\mathcal L}$, we have $${\mathcal E}^*\otimes{\mathcal L}\cong
 {\mathcal E},$$
 Hence (\ref{saumon}) rewrites as
 \begin{eqnarray}
 \label{camille}
 \oplus_iH^0(\tilde R,\tilde g^*(S^{k-3}{\mathcal E}\otimes{\mathcal L})\otimes H_{x_i}^*)
 \rightarrow H^0(\tilde R,\tilde g^*(S^{k-3}{\mathcal E}\otimes {\mathcal E}))
 \end{eqnarray}
 given by the inclusions
 $$H_{x_i}^*\otimes\tilde g^*{\mathcal L}\subset \tilde g^*{\mathcal E}.$$
We want to show that (\ref{sannomnom}) is injective, or that (\ref{camille}) is surjective.
We already noticed that the restriction map
$$H^0(G_2,S^{k-3}{\mathcal E}\otimes {\mathcal E})=S^{k-3}H^0(E)^*\otimes H^0(E)^*$$
$$\rightarrow H^0(R,S^{k-3}{\mathcal E}\otimes {\mathcal E})=
H^0(\tilde R,\tilde g^*(S^{k-3}{\mathcal E}\otimes {\mathcal E}))$$
is surjective.
 On the other hand, consider the $2$-dimensional
 subspace $H_{x_i}=K_{x_i}^{\perp}\subset H^0(E)^*$.
 It is obvious that it parametrizes 
  sections of 
  $${\rm Ker}\,(H^0(P_{x_i},g^*{\mathcal E})\rightarrow H^0(P_{x_i},H_{x_i}))=
  H^0(P_{x_i},g^*{\mathcal L}\otimes H_{x_i}^*).$$
  Hence the surjective map
  $$S^{k-3}H^0(E)^*\otimes H^0(E)^*\rightarrow  H^0(\tilde R,\tilde g^*(S^{k-3}{\mathcal E}
  \otimes{\mathcal E}))$$
  sends 
  $S^{k-3}H^0(E)^*\otimes H_{x_i}$ in   the subspace
  $H^0(\tilde R,\tilde g^*(S^{k-3}{\mathcal E}\otimes{\mathcal L})\otimes H_{x_i}^*).$
 
 Now since the $x_i$'s are generic, the spaces $H_{x_i}$
 generate $H^0(E)^*$, hence
 the $S^{k-3}H^0(E)^*\otimes H_{x_i}$'s generate
 $S^{k-3}H^0(E)^*\otimes H^0(E)^*$. Hence we have shown that
 (\ref{camille}) is surjective.

\cqfd

\section{Appendix}
We consider the Grassmannian $G_2$ of
rank $2$ vector subspaces of a $k+2$-dimensional vector space $V$. Let ${\mathcal L}$
be the line bundle on $G_2$ whose sections give the Pl\"ucker embedding.
If ${\mathcal E}$ is the dual of the tautological subbundle ${\mathcal S}\subset V\otimes
{\mathcal O}_{G_2}$, we have ${\mathcal L}=det\,{\mathcal E}$. The cohomology
groups
$H^p(G_2,{\mathcal L}^{-q}\otimes S^{q'}{\mathcal E})$ are described in the following proposition.
\begin{prop} \label{appen}For $q>0,\,q'>0$, we have 
$$H^p(G_2,{\mathcal L}^{-q}\otimes S^{q'}{\mathcal E})=0\,\,{\rm if}\,\, p\not=k,\,2k.$$
Furthermore, for $p=k$,
we have
$$H^p(G_2,{\mathcal L}^{-q}\otimes S^{q'}{\mathcal E})=0\,\,{\rm if}\,\, -q+q'+1<0,$$ and for
$p=2k$, we have
$$H^p(G_2,{\mathcal L}^{-q}\otimes S^{q'}{\mathcal E})=0\,\,{\rm if }\,\,-q+q'\geq-k-1.$$
\end{prop}
{\bf Proof.}
Let
$$\begin{matrix} &P&\stackrel{g}{\rightarrow }&G_2&\\
&f\downarrow&&&\\
&{\mathbb P}(V)&&&
\end{matrix}
$$
be the incidence variety. $P$ is a ${\mathbb P}^1$-bundle over $G_2$ and
a ${\mathbb P}^k$-bundle over ${\mathbb P}(V)$. Let $H:=f^*{\mathcal O}_{{\mathbb P}(V)}(1)$
and let $L'=g^{*}{\mathcal L}$.
Then
${\mathcal E}=R^0g_*H$ and $S^{q'}{\mathcal E}=R^0g_*(q'H)$.
It follows that we have
$$H^p(G_2,{\mathcal L}^{-q}\otimes S^{q'}{\mathcal E})=
H^p(P,-qL'+q'H).$$
Next the line bundle $L'$ restricts to
${\mathcal O}(1)$ on the fibers of $\pi$. It follows from
this that
$$K_P=-(k+1)L'-2H,$$
and $K_{P/{\mathbb P}(V)}=-(k+1)L'+kH$.

Now since $q>0$ we have $R^lf_*(-qL'+q'H)=0$ for $l<k$ and hence
\begin{eqnarray}
\label{dijan1}
H^p(P,-qL'+q'H)=H^{p-k}({\mathbb P}(V),R^kf_*(L'+q'H)).
\end{eqnarray}
By Serre's duality, we have
$$R^kf_*(-qL'+q'H)=(R^0f_*(qL'-q'H-(k+1)L'+kH))^*$$
$$
=(R^0f_*((q-(k+1))L'+(k-q')H))^*.$$
Now we have 
\begin{eqnarray}
\label{dijan}R^0f_*((q-(k+1))L')=0
\end{eqnarray}
if
$q<k+1$, and 
\begin{eqnarray}
\label{eqomega}
R^0f_*((q-(k+1))L')\cong S^{q-k-1}(\Omega_{{\mathbb P}(V)}(2))
\end{eqnarray}
for $q\geq k+1$.
(The isomorphism (\ref{eqomega}) for $q-k+1=1$ follows
from
the isomorphism
$$H^0(P,L')=H^0(G_2,{\mathcal L})=\bigwedge^2V^*=H^0({\mathbb P}(V),\Omega_{{\mathbb P}(V)}(2))$$
and from the comparison of the kernels of the surjective
evaluation maps
$$H^0(P,L')\rightarrow H^0(f^{-1}(x),L')$$
and
$$H^0({\mathbb P}(V),\Omega_{{\mathbb P}(V)}(2))\rightarrow\Omega_{{\mathbb P}(V)}(2)_x.)$$

Finally we conclude from (\ref{dijan1}), (\ref{dijan}) and
(\ref{eqomega}) that
\begin{enumerate}
\item $H^p(P,-q L'+q'H)=0$ for $p<k$.
\item $H^p(P,-q L'+q'H)=0$ for $q<k+1$.
\item For $p\geq k$, $q\geq k+1$,
$$H^p(P,-qL'+q'H)=H^{p-k}({\mathbb P}(V),S^{q-k-1}(T_{{\mathbb P}(V)}(-2))(q'-k)).$$
\end{enumerate}
To conclude, consider the Euler exact sequence
$$
0\rightarrow {\mathcal O}_{{\mathbb P}(V)}(-1)\rightarrow 
V\otimes{\mathcal O}_{{\mathbb P}(V)}\rightarrow T_{{\mathbb P}(V)}(-1)\rightarrow 0.$$
It induces the exact sequences
$$0\rightarrow  S^{q-k-2}V\otimes{\mathcal O}_{{\mathbb P}(V)}(-q+q')
\rightarrow S^{q-k-1}V\otimes{\mathcal O}_{{\mathbb P}(V)}(-q+q'+1)$$
$$\rightarrow S^{q-k-1}(T_{{\mathbb P}(V)}(-1))(-q+q'+1)\rightarrow 0.$$
Hence we conclude that the space $H^p(G_2,{\mathcal L}^{-q}\otimes S^{q'}{\mathcal E})$
which by the above is also isomorphic to
$$H^{p-k}({\mathbb P}(V),S^{q-k-1}(T_{{\mathbb P}(V)}(-2))(q'-k))$$
$$=H^{p-k}({\mathbb P}(V),S^{q-k-1}(T_{{\mathbb P}(V)}(-1))(-q+q'+1))$$
is equal to $0$ for $p-k\not=0,\,k$ (since $p\leq2k$),
and that :

- for $p-k=0$ it is $0$ if $-q+q'+1<0$;

-  for $p-k=k$ it is $0$
if $-q+q'\geq -k-1$.
\cqfd


\begin{thebibliography}{99}
\bibitem{aprodu} M. Aprodu. On the vanishing of higher syzygies of curves, preprint
2000, to appear in Math. Zeitschrift.
\bibitem{acgh} E. Arbarello, M. Cornalba, Ph. griffiths, J. Harris. 
{\it Geometry of algebraic curves}, Vol. 1, Grundlehren 
der Math. Wissenschaften 267, Springer-Verlag 1985.
\bibitem{} G. Danila.  Sur la cohomologie d'un fibr\'e tautologique
 sur le sch\'ema
de Hilbert d'une surface, Journal of Algebraic Geometry, 10 (2001), no. 2,
p.247-280.
\bibitem{eh} S. Ehbauer. Syzygies of points in projective space and applications,
in {\it Proceedings of the international conference, Ravello 1992},
edited by Orecchia and Chiantini, Walter de Gruyter (1994).
\bibitem{ein} L. Ein. A remark on the syzygies of the generic canonical curve,
J. diff. Geom. 26 (1987), 361-365.
\bibitem{EH} D. Eisenbud, J. Harris. Limit linear series, 
basic theory, Invent. Math. 85 (1986), 337-371.
\bibitem{Gre} M. Green. Koszul cohomology and the geometry of projective varieties, J. Diff. geom. 19,
(1984) 125-171.
\bibitem{GRH} M. Green, R. Lazarsfeld. Special divisors on curves on a $K3$ surface,
Inventiones Math. 89, 357-370 (1987).
\bibitem{georoma} M. Green, R. Lazarsfeld. A simple proof of Petri's theorem on canonical 
curves, in {\it  Geometry today},
 Giornate di Geometria Roma (1984) PM 60, Birkha\"user, 
129-142.
\bibitem{HR} A. Hirschowitz, S. Ramanan. New evidence for Green's conjecture
on syzygies of canonical curves, Ann. Sci. Ecole Norm. sup., IV S\'erie
31(4), 141-152 (1998).
\bibitem{La} R. Lazarsfeld. Brill-Noether-Petri without degenerations, 
J. Diff. Geom. 23, 299-307 (1986).
\bibitem{Mu1} S. Mukai. Biregular classification of Fano 
threefolds 
and Fano manifolds of coindex 3, Proc. Nat. Acad. Sci. U. S. A. 86 
(1989) 3000-3002.
\bibitem{Mu2} S. Mukai. Symplectic structure on the moduli space
of sheaves on an abelian or $K3$ surface, Invent. Math. 77, 101-116 (1984).
\bibitem{Sch1} F.-O. Schreyer.
 Syzygies of canonical curves and special linear series, Math. Ann. 275 (1986), 
105-137.
\bibitem{Sch2}  F.-O. Schreyer. A standard basis approach to
 the syzygies of canonical curves, J. Reine und Angew. Math. 421 (1991), 83-123.
\bibitem{Teix} M. Teixidor.  Green's conjecture for 
the generic canonical curve, preprint 1999.
\bibitem{te} M. Teixidor.... Duke Math. Journal (2002).
\bibitem{voi1} C. Voisin. Sur l'application de Wahl des courbes satisfaisant la condition 
de Brill-Noether-Petri, 
Acta Mathematica , vol. 168 (1992) 
249-272.
\bibitem{voi2} C. Voisin. Courbes t\'etragonales et cohomologie de Koszul, \newline
J. Reine Angew. 
Math. 387 (1988).
\bibitem{voi3} C. Voisin. D\'eformation des syzygies et th\'eorie de Brill-Noether,\newline
 Proc. 
London Math. Soc. 3 ,vol. 67 (1993) 493-515.

\end{thebibliography}
\end{document}